\documentclass{article}
\usepackage[utf8]{inputenc}
\usepackage[margin=1in]{geometry}
\usepackage{graphicx}
\usepackage{amsmath}
\usepackage{subfigure}
\usepackage{algorithm}
\usepackage{algorithmic}

\usepackage[square,numbers]{natbib}
\usepackage{caption}
\usepackage[version=4]{mhchem}
\usepackage[dvipsnames]{xcolor}
\usepackage{float}
\usepackage{amsthm}
\usepackage{halloweenmath}
\usepackage[affil-it]{authblk}
\usepackage{amssymb}

\newtheorem{proposition}{Proposition}

\theoremstyle{definition}

\newtheorem{remark}{Remark}

\title{Hunting $\varepsilon$: The origin and validity of quasi-steady-state 
reductions in enzyme kinetics}
\author[1,X]{Justin Eilertsen}
\author[1,Y]{Malgorzata Anna Tyczynska}
\author[1,2,Z]{Santiago Schnell}
\affil[1]{Department of Molecular \& Integrative Physiology, 
University of Michigan Medical School, Ann Arbor, Michigan, USA}
\affil[2]{Department of Computational Medicine \& Bioinformatics, 
University of Michigan Medical School, Ann Arbor, Michigan, USA}
\affil[X]{Current affiliation: Mathematical Reviews, American Mathematical Society, 
416 $4th$ Street, Ann Arbor, MI 48103}
\affil[Z]{Current affiliation: Integrated Mathematical Oncology Department,
H. Lee Moffitt Cancer Center and Research Institute,
12902 USF Magnolia Drive, SRB-4, Tampa, FL 33612}
\affil[Z]{Corresponding author. Current affiliation: Department of Biological Science and 
Department of Applied and Computational Mathematics and Statistics, University of 
Notre Dame, Notre Dame, IN 46556, USA. E-mail: santiago.schnell@nd.edu}

\date{}

\begin{document}

\maketitle
\begin{abstract}
The estimation of the kinetic parameters that regulate the speed of enzyme catalyzed
reactions requires the careful design of experiments under a constrained set of 
conditions. Many estimates reported in the literature incorporate protocols that leverage
simplified mathematical models of the reaction's time course known as quasi-steady-state
reductions. Such reductions often --- but not always --- emerge as the result of a singular
perturbation scenario. However, the utilization of the singular perturbation reduction method 
requires knowledge of a dimensionless parameter, ``$\varepsilon$", that is proportional 
to the ratio of the reaction's fast and slow timescales. To date, no such ratio has been
determined for the intermolecular autocatalytic zymogen activation reaction, which means 
it remains open as to when the experimental protocols described in the literature are 
even capable of generating accurate estimates of pertinent kinetic parameters. Using 
techniques from differential equations, Fenichel theory, and center manifold theory, 
we derive the appropriate ``$\varepsilon$" whose magnitude regulates the validity of 
the quasi-steady-state reduction employed in the reported experimental procedures. 
Although the model equations are two-dimensional, the fast/slow dynamics are rich. The 
phase plane exhibits a dynamic transcritical bifurcation point in a particular singular 
limit. The existence of such a bifurcation is relevant, because the critical manifold 
losses normal hyperbolicity and classical Fenichel theory is inapplicable. We show 
that while there exists a faux Canard that passes directly through this bifurcation 
point, trajectories emanating from experimental initial conditions are actually bounded 
away from the bifurcation point by an asymptotic distance that is proportional to 
the square root of the linearized system's eigenvalue ratio. Furthermore, we show 
that in some cases chemical reversibility can be interpreted dynamically as an 
imperfection, since the presence of reversibility can destroy the bifurcation 
structure present in the singular limit. By extension, some of these features are 
also present in the phase--plane dynamics of the famous Michaelis--Menten reaction 
mechanism. Finally, we show that the reduction method by which QSS reductions are 
justified can depend on the path taken in parameter space. Specifically, we show 
that the standard quasi-steady-state reduction for this reaction is justifiable by 
center manifold theory in one limit, and via Fenichel theory in a different limit. 
\end{abstract}
\newpage 

\section{Introduction}
The model equations of biochemical systems that exhibit disparate timescales can often 
be systematically reduced. Classical Fenichel theory \cite{Fenichel1979} rigorously 
justifies model reduction via projection onto a slow, normally hyperbolic invariant 
manifold, thereby reducing the number of equations needed to capture relevant 
biophysical chemical events that occur over slow timescales. The extension of singular
perturbation theory to systems whose critical manifolds lose normal hyperbolicity has 
led to our understanding of nonlinear phenomena in biological systems~\cite{BERTRAM2017},
including relaxation oscillations~\cite{Krupa2001}, mixed--mode 
oscillations~\cite{MMO,MMO1,MM03}, and cardiac arrhythmias~\cite{Kimrey2020}. However,
there is another, less emphasized --- but equally important --- avenue for singular 
perturbation theory in interdisciplinary science and applied mathematics: 
\textit{experimental design}. In chemical kinetics, reduced models that materialize 
from the application of Geometric Singular Perturbation Theory (GSPT) are called
\textit{quasi-steady-state approximations} (QSSAs), or QSS \textit{reductions}. In 
enzyme kinetics, QSS reductions are a predominant tool in functional characterization 
and drug design: kinetic parameters are estimated by fitting experimental time course 
data to QSS reductions via nonlinear regression 
analysis \cite{Schnell2003Estimates,Stroberg:2016,choiandKim2017}. In turn, 
the estimated kinetic parameters provide valuable information concerning the viability 
and efficacy of an enzyme, substrate, inhibitor, or activator in a reaction mechanism.
Accordingly, GSPT plays a central role in metabolic engineering and drug design.

While the kinetic parameters of a particular reaction mechanism will be more or less
invariant with respect to certain thermodynamic constraints, the QSSAs employed in 
experiments are approximations and, as such, will generally be \textit{invalid} 
unless special conditions are upheld. Consequently, if we are to estimate intrinsic kinetic
parameters by fitting data to QSSAs, experiments must be prepared in a way that is 
compatible with the validity of the QSS reduction \cite{Schnell2003Estimates}. 

Since the application of GSPT relies on the existence fast and slow timescales, 
the experimental constraint that secures the legitimacy of a QSS reduction has a 
direct mathematical translation: The accuracy of a QSS reduction depends on the 
ratio of the fast and slow timescales, $\varepsilon$, which must be small to 
ensure the validity of the QSS reduction. Experiments must therefore be prepared 
so that $0< \varepsilon \ll 1$. But, herein lies the two-fold caveat: how do we 
know that a QSS reduction is the result of a singular perturbation and, if it is, 
then what \textit{is} $\varepsilon$? This question might be superfluous to the 
theoretician whose job is develop and expand the theory of singular
perturbations. In contrast, this is \textit{the} vital question posed to the
interdisciplinary scientist whose objective is to use QSS reductions to extract 
quantitative information from experimental data.

To articulate the challenges of this problem, and profile how mathematics plays a 
crucial role, we have chosen to analyze the intermolecular autocatalytic zymogen 
activation (IAZA) mechanism, represented schematically by
\begin{equation}\label{z1}
    \ce{$Z$ + $E$ <=>[$k_1$][$k_{-1}$] $C$ ->[$k_2$] $2E$ + $W$},
\end{equation}
where $Z$ is a zymogen, $E$ is an active enzyme, $W$ is peptide, and $k_1,k_2$ 
and $k_{-1}$ are deterministic rate constants. The chemical mechanism (\ref{z1}) 
is essentially the autocatalytic ``cousin" of the famous Michaelis--Menten reaction
mechanism. As enzyme precursors (proenzymes), zymogens have numerous biochemical 
functions. For example, they play a critical role in protein digestion by converting 
pepsin to pepsinogen~\citep{al1972kinetics,sanny1975conversion}, 
and trypsin to trypsinogen~\citep{kassell1973zymogens,khan1998molecular,thrower2006zymogen}.
However, the real motivation behind our selection of the IAZA mechanism for analysis 
emanates from the fact that there is a generous body of experimental literature 
reporting on parameter estimates extracted from fitting data to reduced 
models~\cite{Wu2001,fuentes2005kinetics,Garcia-Moreno1991}. However, the precise, 
theoretical qualifiers that warrant the validity of the QSS reductions employed in 
the experimental literature is open and unresolved. Unfortunately, this means that 
the accuracy of estimated parameters reported in the literature depends on, at least 
in part, whether or not the experiments where prepared in a way that justifies the 
use of the QSSA employed in the regression procedure~\cite{hanson2008RSA}.


\section{The IAZA reaction mechanism}\label{sec:phase-plane}
In this section we introduce the model equations of the IAZA reaction mechanism. 
We also give a brief review of the reduced models employed in the experimental 
literature, and discuss two \textit{classical} QSSAs commonly attributed to GSPT 
and slow manifold projection.
 
\subsection{The mass action model equations}
To derive a deterministic (and finite-dimensional) model of the IAZA reaction
mechanism~(\ref{z1}), we let $z:=z(t)$, $e:=e(t)$ and $c:=c(t)$ denote the 
concentrations of $Z$, $E$, and $C$ respectively, and apply the law of mass action 
to (\ref{z1}). This yields the following set of nonlinear ordinary differential 
equations
\begin{subequations}
\begin{align}
\dot{z} &=-k_1 ez +k_{-1}c,\label{eq:ma1}\\
\dot{c} &=k_1 ez - (k_{-1}+k_2)c,\label{eq:ma2}\\
\dot{e} &=-k_1 ez + (k_{-1}+2k_2)c,\label{eq:ma3}\\
\dot{w} &=k_2c, \label{eq:ma4}
\end{align}
\end{subequations}
with ``$\dot{\phantom{x}}$" denoting differentiation with respect to time. The 
typical initial conditions used in experiments are
\begin{equation} \label{z3} 
    \begin{matrix}
    z(0) = z_0, & & e(0) = e_0, & & c(0) = 0, & & w(0)=0,
    \end{matrix}
\end{equation}
and, unless otherwise stated, we will assume \eqref{z3} holds in the analysis 
that follows. 

The structure of the model equations implies the existence of conservation laws.  
Note that equations (\ref{eq:ma1})--(\ref{eq:ma3}) obey
\begin{align}
\dot{z} + \dot{c}+ \dot{w} =0,\label{eq:con1}\\
\dot{e} + \dot{c} - \dot{w} =0,\label{eq:con2}
\end{align}
and their sum yields an equation that is independent of $w$: 
\begin{equation}\label{eq:con}
\dot{z} + \dot{e} + 2\dot{c}=0.
\end{equation}
Integrating (\ref{eq:con}) with respect to (\ref{z3}) yields the conservation law 
$z+e+2c=E_T$, where $E_T$ is sum of initial concentrations of zymogen and enzyme: 
$E_T=z_0+e_0$. Substituting $e=E_T-z-2c$ into 
(\ref{eq:ma1})--(\ref{eq:ma2}) yields
\begin{subequations}
\begin{align} 
   \dot{z} &= -k_1(E_T - z)z + (k_{-1} + 2k_1z)c,\label{z5a}\\
    \dot{c} &= k_1(E_T -z)z - (k_{-1}+k_2 + 2k_1z)c ,  \label{z5b}\\
    \dot{w} &= k_2 c.\label{z5c}
\end{align}
    \end{subequations} 
Summing together (\ref{z5a})--(\ref{z5c}) 
and integrating yields the additional conservation law
\begin{equation}\label{eq:Zycon}
z_0 = z+c+w,
\end{equation}
and thus we can express the mass actions equations in terms of two variables: 
$(z,c), (w,z)$ or $(w,c)$.

In what follows, we will occasionally refer to the parameters $K_M,K_S$ and $K$, 
which are comprised of rate constants. The parameters $K_M,K_S$, and $K$ are:
\begin{equation}
    K_M:=\cfrac{k_{-1}+k_2}{k_1}, \qquad K_S:=\cfrac{k_{-1}}{k_1}, \qquad K:=\cfrac{k_2}{k_1}.
\end{equation}
Note that $K_M=K_S+K$.

\subsection{Phase plane geometry}
\begin{figure}[htb!]
  \centering
    \includegraphics[width=7cm]{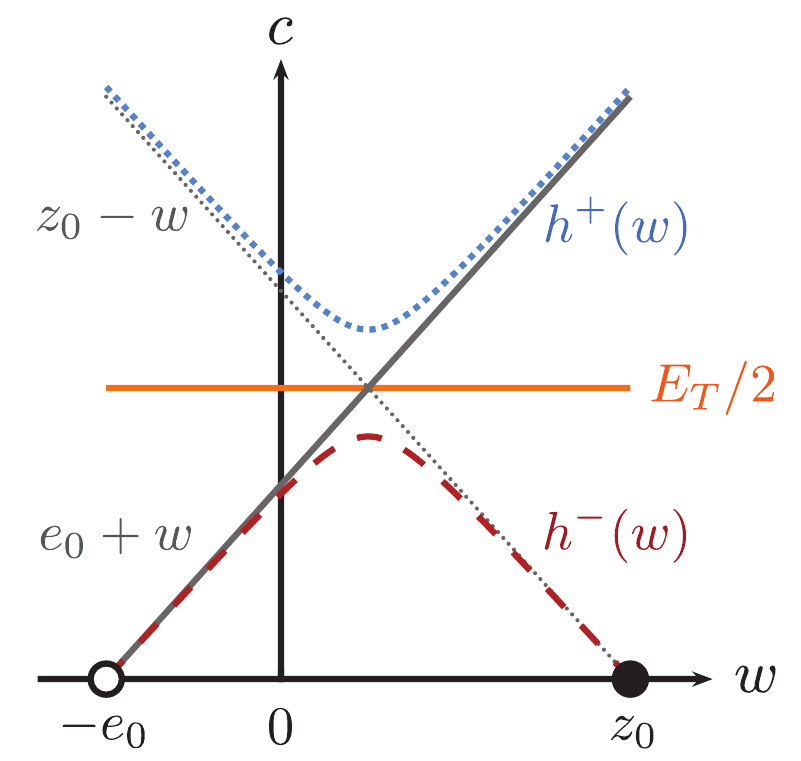}
  \caption{\textbf{The phase plane geometry associated with the IAZA reaction 
  mechanism~(\ref{z1}) contains two equilibrium points.} The dashed red curve is 
  the $c$-nullcline given by $c=h^-(w)$, and the dotted blue curve corresponds to
  $c=h^+(w)$. The horizontal orange line is the constant $c=E_T/2$. The gray 
  diagonal dotted line is $c=z_0-w$, and the gray diagonal solid line is $c=w+e_0$. 
  The $c$-nullcline coincides with the $w$-axis at two equilibrium points. The black 
  circle is the attracting equilibrium point $x^{(1)}$, located at $(z_0,0)$ and the 
  saddle equilibrium $x^{(2)}$ is located at $(-e_0,0)$ and marked by a white
  circle.}\label{fig0}
\end{figure}

The mass action equations in $(w,c)$ coordinates are formulated from the conservation 
law (\ref{eq:Zycon}):
\begin{subequations}
\begin{align}
    \dot{c} &= k_1(e_0+w-c)(z_0-c-w)-(k_{-1}+k_2)c,\label{eq:cdotMA}\\
    \dot{w} &= k_2c\label{eq:wdotMA}.
\end{align}
\end{subequations}
Note that there are three nullclines in the $(w,c)$ phase--plane, $c=0,$ and $c=h^{\pm}(w)$:
\begin{subequations}\label{cnulls}
\begin{align}
    c=h^-(w) &:= \cfrac{1}{2}(K_M+E_T)-\cfrac{1}{2}\sqrt{(K_M+E_T)^2-4(z_0-w)(e_0 + w)}, \label{eq_h-} \\
    c=h^+(w) &:= \cfrac{1}{2}(K_M+E_T)+\cfrac{1}{2}\sqrt{(K_M+E_T)^2-4(z_0-w)(e_0 + w)}. \label{eq_h+}
\end{align}
\end{subequations}
The nullcline $c=h^-(w)$ intersects with $c=0$ at the points
\begin{equation}
    x^{(1)} = (z_0,0), \quad x^{(2)} = (-e_0,0).
\end{equation}
Linearization reveals that $x^{(1)}$ is an attracting equilibrium point and $x^{(2)}$ is 
a saddle point\footnote{One can speculate about the possible existence of a heteroclinic 
connection between $x^{(1)}$ and $x^{(2)}$. However, we not not pursue this in any 
detail here, as it is not critical to our analysis.} (see {\sc{Figure}}~\ref{fig0}). 

\begin{figure}[htb!]
  \centering
    \includegraphics[width=14cm]{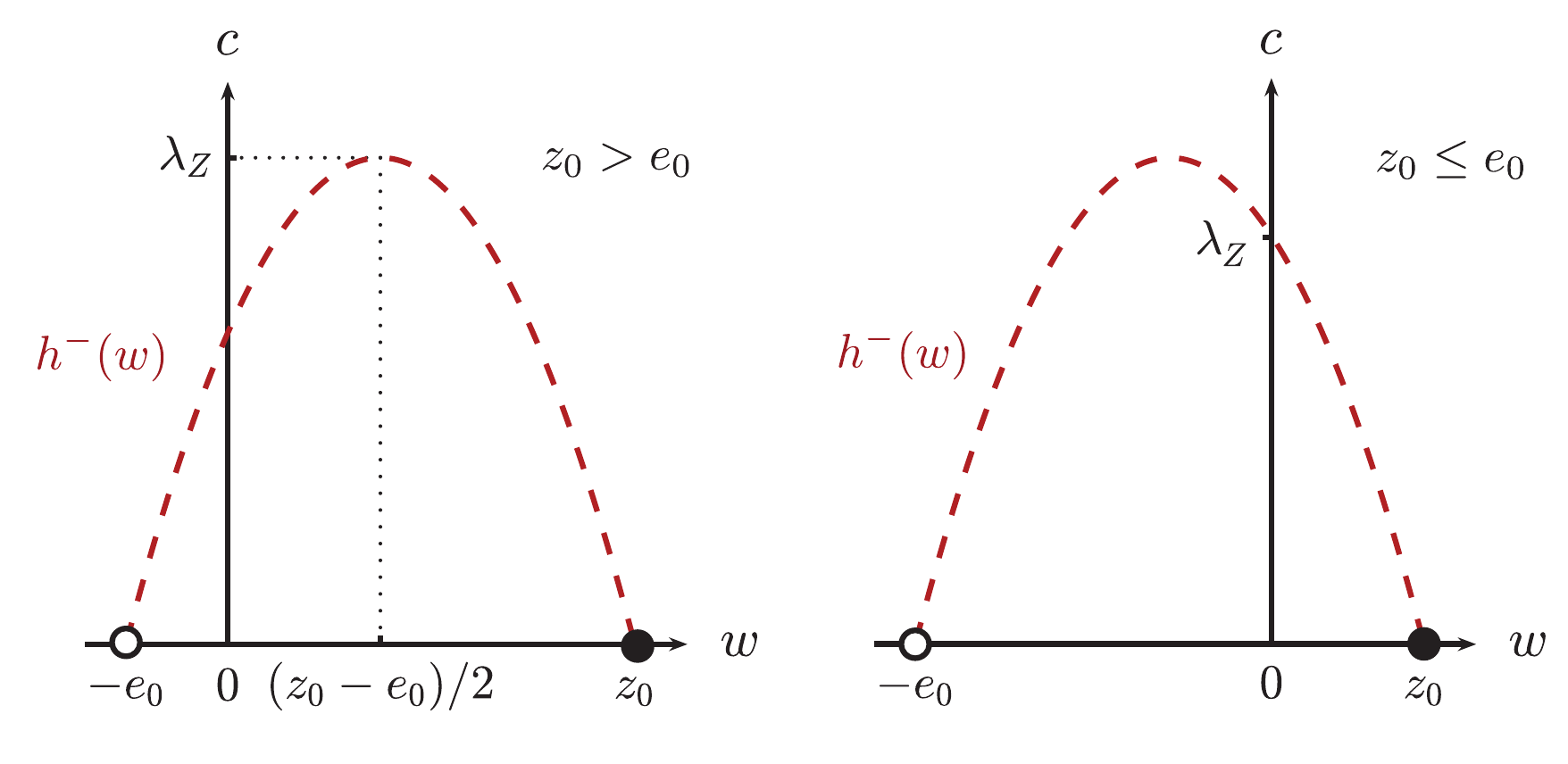}
  \caption{\textbf{The supremum of $\mathbf{c}$ depends on the initial ratio of enzyme 
  to zymogen concentration in the IAZA reaction   mechanism~(\ref{z1}).} In both panels, 
  the dashed, red curve is the $c$-nullcline 
  given by $c=h^-(w)$. {\sc Left}: This is an illustration of the phase plane when 
  $e_0 < z_0$. In this case, the apex of the nullcline given by $c=h^-(w)$ occurs at 
  $w_T=(z_0-e_0)/2$, which is positive. Consequently, the complex concentration is 
  bounded by $c \leq \lambda_Z = h^-(w_T)$. {\sc Right}: In this illustration, the 
  initial enzyme concentration is greater than the initial zymogen concentration, 
  i.e., $z_0 < e_0$. The apex of the $c$-nullcline now lies in quadrant II, and 
  corresponds to a negative and unphysical $w_T$. Consequently, $c$ is bounded 
  above by $h^-(0)$ when $z_0 < e_0$. }\label{fig0A}
\end{figure}

When typical initial conditions for experiments (\ref{z3}) are imposed, the nullcline 
given by $c=h^-(w)$ provides an upper bound on the complex concentration, $c$. There 
are two cases we must work out. The first case we consider is $e_0 \leq z_0$. 
In this scenario, the complex concentration has a maximum given by the apex of the 
curve $c=h^-(w)$. Differentiating this expression with respect to $w$ yields a critical 
point at $w=(z_0-e_0)/2\equiv w_T$, and thus $ \max c \leq \lambda_Z:= h^-(w_T)$ 
(see {\sc Figure}~\ref{fig0A}, {{\sc left}} panel).

In contrast, if $z_0 < e_0$, then the critical point at $w_T$ is negative, and therefore 
nonphysical. Thus, when experimental initial conditions are prescribed and $z_0<e_0$, 
the supremum of $c$ is given by $h^-(0)$, since $c=h^-(w)$ is a monotonically decreasing 
function on the interval $0\leq w\leq z_0$ (see {\sc Figure}~\ref{fig0A}, {{\sc right}}
panel).
From this point forward, we will simply refer to the maximum of $c$ as $\lambda_Z$, 
and take this to be $h^-(w_T)$ if $e_0 \leq z_0$, or $h^-(0)$ if $z_0 < e_0$.

\begin{figure}[htb!]
  \centering
    \includegraphics[width=8cm]{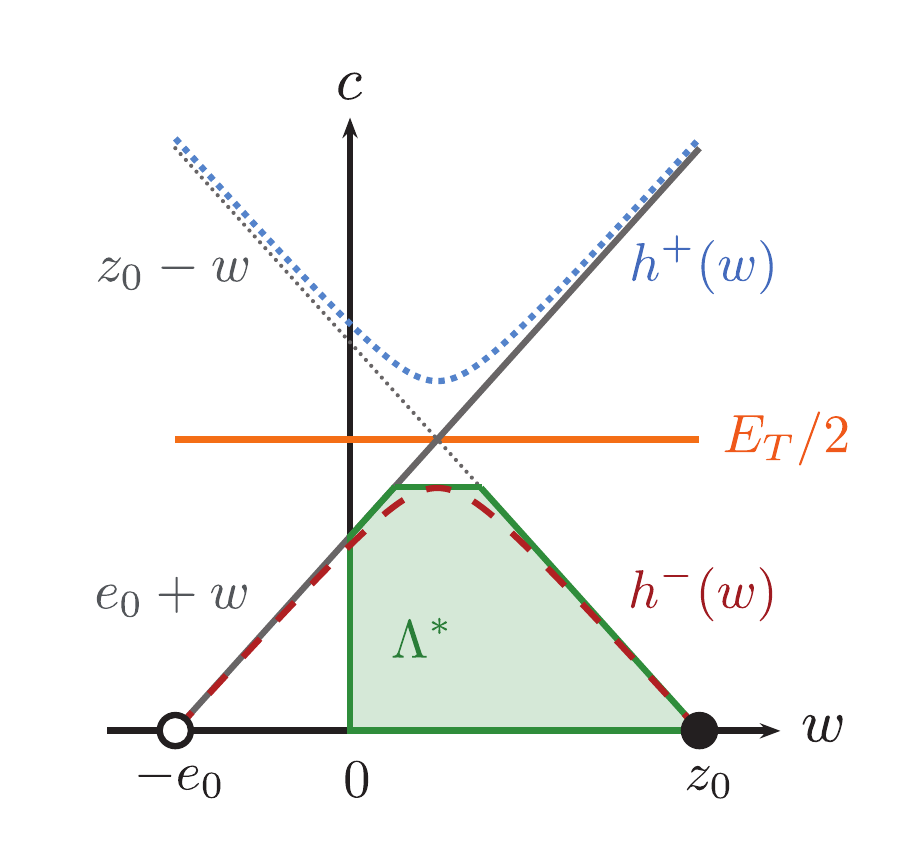}
  \caption{\textbf{Due to the presence of conservation laws and attractors, the 
  phase plane contains forward invariant sets in the IAZA reaction 
  mechanism~(\ref{z1})}. The set $\Lambda^*$ is filled with green shade and outlined 
  by the green lines. Any trajectory that starts on or in   the $\Lambda^*$ stays 
  on or in it for all $t>0$.}\label{invariant}
\end{figure}

Note that due to conservation, any trajectory that starts inside or on the region
``$\Lambda$" bounded by the curves
\begin{equation}\label{eq:curves}
\Lambda := \{(w,c)\in \mathbb{R}_{\geq 0}^2 :  c \leq E_T/2\} \cap \{(w,c)\in \mathbb{R}_{\geq 0}^2 : c \leq e_0+w\} \cap \{(w,c)\in \mathbb{R}_{\geq 0}^2 :c \leq z_0-w\}
\end{equation}
must stay on the boundary or inside $\Lambda$ for all positive time. Thus, $\Lambda$ is positively 
invariant. Moreover, we can be even more restrictive, and define $\Lambda^*$ to be
\begin{equation}
    \Lambda^*:=\Lambda \cap \{(c,w)\in \mathbb{R}^2_{\geq 0}: c \leq \lambda_Z\},
\end{equation} 
which is also positively invariant (see {\sc Figure}~\ref{invariant} for an 
illustration of $\Lambda^*$).

\subsection{Strategies for model reduction: A brief literature review}
Although the complete system (\ref{eq:ma1})--(\ref{eq:ma4}) is expressible in 
terms of only two variables, there is a two-part question that we need to answer: 
Is it possible to reduce the system (\ref{z5a})--(\ref{z5c}) further and, if so, 
how favorable is such a reduction with respect to a range of parameters? The 
answer to this question is critical for the derivation of approximate mathematical 
models that are more favorable (well-suited) for estimating the pertinent kinetic 
parameters (i.e., rate constants). Several authors attempted to address this very 
question. An early attempt to construct a reduced model for the IAZA reaction was 
made by \citet{Garcia-Moreno1991}. These authors generated what is known as a 
pseudo-first-order (PFO) approximation, and assumed that if $e_0 \ll z_0$, 
then $z\approx z_0$ during the initial phase of the reaction. Replacing $z$ 
with $z_0$ in (\ref{z5a})--(\ref{z5c}) results in a set of linear equations that 
can be easily be solved. Unfortunately, there is a downside to the PFO approach: 
the approximation is only valid (if at all) at the onset of the reaction, and 
therefore its accuracy is limited to very short (brief) timescales. Nevertheless, 
and despite the limitations of the linear approximation that arises from the PFO 
procedure, the analysis conducted in \citet{Garcia-Moreno1991} was essentially 
the first study in which considerable progress was made towards our understanding 
of autocatalytic zymogen activation reaction kinetics, since  prior to the work 
of \citet{Garcia-Moreno1991}, mostly ``over-simplified" reaction mechanisms 
describing zymogen activation were studied \citep{tans1983,tans1987}. Later 
investigations \citep{fuentes2005kinetics,varon2006competitive,varon2006kineticJoTB} 
added complexity to the IAZA reaction mechanism (\ref{z1}), but the kinetic 
analyses conducted in these studies employed the basic PFO methodology 
of \citet{Garcia-Moreno1991}. 

A slightly different and more rigorous approach was developed by \citet{Wu2001}. 
Instead of constructing a reduced model by assuming PFO kinetics, \citet{Wu2001} 
derived a simplified model by employing a QSSA that presumes the complex concentration, 
$c$, changes very slowly ($\dot{c}\approx 0)$ for the majority of the reaction. 
By leveraging their reduced model based on the QSSA, \citet{Wu2001} developed 
an experimental protocol and found that their estimated parameter values were 
close to those obtained in an earlier study by \citet{Garcia-Moreno1991}. Additional 
studies utilized the protocol of \cite{Wu2001} to estimate kinetic parameters for 
ligand-induced autocatalytic reactions \citep{WangLigand}, and autocatalytic 
reactions that contain competitive inhibitors \citep{Wang2004Inhibitor}. However, 
despite the promising results that emerged from the work of \citet{Wu2001}, the 
authors did not establish a suitable qualifier that guarantees the validity of 
their QSSA. Thus, it is impossible to determine the appropriate bounds within 
which their protocol for estimating kinetic parameters is reliable. 
 
\subsection{The \textit{classical} QSS reductions of the IAZA reaction}
As we mentioned earlier, QSS reductions in enzyme kinetics are based on 
the assumption that the reaction will undergo an initial transient phase in which 
the complex concentration accumulates rapidly, followed by a QSS phase where the 
rate of change in complex concentration nearly vanishes. The assumption that 
$\dot{c}$ is nearly zero provides a heuristic avenue for model reduction. Solving 
the equation $k_1(E_T-z)z -(k_{-1}+k_2+2k_1z)c=0$ for $c$ yields
\begin{equation}\label{DAEC}
    c=\cfrac{(E_T-z)z}{K_M+2z},
\end{equation}
and insertion of (\ref{DAEC}) into (\ref{z5a}) generates the \textit{classical} 
QSS reduction for $z$:
\begin{equation}\label{zclass}
    \dot{z} = -\cfrac{k_2}{K_M+2z}(E_T-z)z.  
\end{equation}
Similarly, in $(c,w)$ coordinates, we take $c=h^-(w)$, which yields the another
type of QSS reduction\footnote{In enzyme kinetics, this particular type of  
reduction is referred to a the \textit{total} quasi-steady-state approximation 
(tQSSA). See, references \cite{BORGHANS1996,tzafriri2003,2002SchnellMaini,Bersani2015} 
for tQSSA analysis in enzyme kinetics.} for $w$, which is equivalent to the 
reduction employed by \citet{Wu2001}:
\begin{equation}\label{wclass}
    \dot{w}=k_2h^-(w).
\end{equation}

Since the QSS phase of the reaction is valid on the slow timescale, the reductions
(\ref{zclass}) and (\ref{wclass}) are generally assumed to hold for the majority 
of the reaction, and kinetic parameters are estimated by fitting (\ref{zclass}) 
or (\ref{wclass}) to experimental time course data; see \citet{Wu2001} as 
an example that deals explicitly with the IAZA mechanism. For a more general 
discussion of mathematical estimation strategies, we invite the reader to 
consult~\cite{Stroberg:2016,choiandKim2017}.

As previously stated, the classical reductions (\ref{zclass}) and (\ref{wclass}) 
are generally valid when there is an initial transient in which $c$ accumulates 
rapidly, followed by QSS phase in which $\dot{c}$ is very small. The fast and 
slow timescales measure the approximate duration of each phase and, while 
calculating classical QSS reductions is rather straightforward, difficulty 
generally emerges in the quest to justify (\ref{zclass}) and (\ref{wclass}). 
Ultimately, justification of (\ref{zclass}) or (\ref{wclass}) is the crux of 
the problem, especially if the objective is to prepare an experiment so 
that (\ref{zclass}) or (\ref{wclass}) is valid. The overwhelming tendency is to 
attribute the validity of an accurate QSS reduction to a singular perturbation 
scenario via scaling analysis. The motivation is that the qualifier (i.e., the 
timescale ratio
``$\varepsilon$") whose smallness justifies the QSS reduction usually emerges 
quite naturally from the scaled equations \cite{Heineken1967,segel1988}. However, 
the difficulty is that \textit{good} scaling relies on accurate timescale 
estimates \cite{segel1988,Segel1989}. Reliable and tractable timescale estimates can be 
difficult to obtain\footnote{Eigenvalues, even in two-dimensional systems, can 
often assume a form that is rather cumbersome.} and, even if the scaling analysis 
suggests a singular 
perturbation scenario in which the theorems of \citet{tikhonov1952} 
and \citet{Fenichel1979} appear to apply, the validity of a QSSA may not actually 
be the result of a singular perturbation scenario. As an example, we invite the
reader to consult \citep{SCHAUER1979A,Noethen:2009:QSS,GOEKE2017,Eilertsen2021}, 
as well as \cite{Goeke2012}, Section~\ref{DTFPV}, for a discussion on scaling analysis. 
This begs the question: Is there a diagnostic methodology capable of determining 
when QSS reduction is the result of a singular perturbation? This answer to this 
question is yes, and the specific methodology that helps to delineate the 
applicability of Fenichel theory, referred to a Tikhonov--Fenichel 
Parameter Value (TFPV) theory, was recently developed by \citet{GOEKE2017}. 
We give a brief overview of TFPV theory in Section~\ref{GSPT}.

\section{Application of the geometric singular perturbation theory and 
Tikhonov--Fenichel parameter values}\label{GSPT}
Our motivation is to generate QSS reductions for the IAZA reaction mechanism 
from GSPT. While GSPT is certainly not the only reduction method one can employ 
to derive QSS reductions, it is favorable since it ensures the existence disparate 
fast and slow timescales. From an experimental vantage point the presence of fast 
and slow timescales is advantageous: the transient window that precedes the QSS 
phase will be short in comparison to the duration of the reaction, which means 
the QSSA will be valid for practically the entire time course, thereby ensuring a
large window from which to extract functional data from experiments.

In this section, we review the basics of GSPT for systems that are not in standard 
fast/slow form. While these results are well-established, they are not as prevalent 
in the literature. Our outline primarily follows~\cite{Martin2020}, Chapter 3. For 
further details, we invite the reader to consult~\cite{Fenichel1979,Goeke2012,Goeke2014}. 
We also give an overview of the more recent TFPV theory developed by 
\citet{GOEKE2017}, which establishes criteria for the applicability of GSPT for 
enzymatic reactions. We conclude with a brief discussion on the missing piece of 
TFPV theory, and the open question that needs resolution.

\subsection{Geometric Singular Perturbation Theory}
To calculate QSS reductions we will employ GSPT. For a perturbed dynamical system 
of the form
\begin{equation}
    \dot{x} = H(x) = h(x) + \varepsilon G(x,\varepsilon), \qquad x\in\mathbb{R}^n, \quad h:\mathbb{R}^n \mapsto \mathbb{R}^n, \quad G:\mathbb{R}^n \times \mathbb{R} \mapsto \mathbb{R}^n, \quad 0 < \varepsilon \ll 1,
\end{equation}
the direct application of GSPT requires the existence of a normally hyperbolic, 
connected, and differentiable manifold of non-isolated stationary points, $M$, in 
the singular limit that corresponds to $\varepsilon =0$:\footnote{In (\ref{CM}),
$\mathbb{O}^n$ denotes an $n$-dimensional column vector with all components 
identically zero, and $Dh(x)$ is the Jacobian of $h$.}
\begin{equation}\label{CM}
    M:=\{ x \in \mathbb{R}^n : h(x)=\mathbb{O}^n\}, \qquad \dim M = k, \quad k< n.
\end{equation}
The manifold $M$ is called the critical manifold. For $x\in M$, the $k$ 
zero-eigenvalues of $Dh(x)$ are referred to as \textit{trivial} eigenvalues, 
and the remaining $n-k$ eigenvalues are called the \textit{non-trivial} 
eigenvalues. Unless otherwise stated, we will assume that the real parts of the nontrivial eigenvalues are strictly less than zero, so that $M$ is attracting. Furthermore, the tangent space of $M$ at $x$, $T_xM$, is given 
by the kernel of $Dh(x)$:
\begin{equation}
    T_xM=\ker Dh(x).
\end{equation}

Fenichel theory applies exclusively to compact subsets, $M_0$, of $M$: 
$M_0\subset M$. Normal hyperbolicity ensures the existence of the splitting,
\begin{equation}\label{split}
    \mathbb{R}^n = T_xM_0 \oplus N_x, \qquad \forall x \in M_0,
\end{equation}
where the subspace $N_x$ in (\ref{split}) is identically the range of $Dh(x)$ 
and complementary to $T_xM$. This 
structure implies the existence of a projection operator, $\Pi^M$,
\begin{equation}
    \Pi^M: \mathbb{R}^n \mapsto T_x M_0,
\end{equation}
and the leading order QSS reduction is obtained by projecting the perturbation,
$\varepsilon G(x,0)$, onto the tangent space of $M_0$:\footnote{We invite the
reader to consult \cite{Goluskin} for a projection method that appears to be 
accurate in special cases.}
\begin{equation}
    \dot{x} = \Pi^{M} G(x,0)|_{x\in M_0}.
\end{equation}

To construct $\Pi^M$, one employs the factorization $h(x) = N(x)f(x)$, where the 
zero level set of $f(x)$ corresponds to $M$, the columns of $N(x)$ for a basis 
for the range of $Dh(x)$, and the row vectors of $Df(x)$ form a basis for the 
orthogonal complement of $\ker Dh(x)$:
\begin{equation}\label{spaces}
\text{column space}\; N= \text{range} \; Dh(x), \qquad \text{column space} \;Df^T= \ker^{\perp} Dh(x),\qquad \forall x \in M_0.
\end{equation}
It follows directly from (\ref{spaces}) that
\begin{equation}
    \Pi^M:=\mathbb{I}^{n\times n} - N (Df N)^{-1} Df,
\end{equation}
where $\mathbb{I}^{n\times n}$ denotes the $n \times n$ identity matrix.

As a final remark, we note that the application of GSPT is usually carried out 
when the perturbed dynamical system is in \textit{standard form}. For a fast/slow 
system in $\mathbb{R}^2$, the standard form is:
\begin{equation}
    \begin{bmatrix} \dot{x} \\ \dot{y} \end{bmatrix} = \begin{bmatrix} 0\\ h(x,y)\end{bmatrix} + \varepsilon \begin{bmatrix} g_1(x,y,\varepsilon) \\ g_2(x,y,\varepsilon) \end{bmatrix}.
\end{equation}
In this 
special case, we have
\begin{equation}
    Nf= \begin{bmatrix} 0 \\ 1 \end{bmatrix}h(x,y),\quad Df \equiv [ D_x h\;\;D_y h ],\quad \Pi^M:=\begin{bmatrix} 1 & 0 \\ -(D_yh)^{-1}(D_xh) & 0\end{bmatrix},
\end{equation}
and the QSS reduction is given by
\begin{equation}
\begin{bmatrix} \dot{x} \\ \dot{y} \end{bmatrix} = \begin{bmatrix} g_1(x,y,0)\\ -(D_yh)^{-1}(D_xh) g_1(x,y,0)\end{bmatrix}.
\end{equation}
Normal hyperbolicity, $D_yh \neq 0$, ensures --- by the Implicit Function Theorem --- 
that $M_0$ is locally expressible as a graph over 
the slow variable, $x$: $y=X(x)$ with $h(x,X(x))=0$, and the QSS reduction 
in standard form is
\begin{equation}
    \dot{x}= g_1(x,X(x),0).
\end{equation}

Often a coordinate transformation can be utilized to recover the standard form 
(see \cite{Martin2020}, Section 3.7, Lemma 3.5), and  well-established results 
in enzyme kinetics have made use of the standard form. However, the standard 
form in the applied literature is generally recovered \textit{not} via coordinate
transformation but via scaling and non-dimensionalization of the mass-action 
system~\cite{Segel1989}. The problem with any scaling procedure is that 
dimensionless variables are not unique and, despite the narrative established 
by several classic papers on singular perturbation analysis in enzyme 
kinetics~\cite{Heineken1967,segel1988}, the standard form appears to be the 
exception, not the rule. 

\subsection{Tikhonov-Fenichel Parameter Values}
As noted in the previous subsection, singularly perturbed mass action systems 
in enzyme kinetics are not generally in standard form. Furthermore, valid QSS 
reductions are not necessarily the result of a singular perturbation
scenario~\cite{Eilertsen2021,GOEKE2017,SCHAUER1979A,Noethen:2009:QSS}. This is 
not at all obvious. \citet{GOEKEpara,GOEKE2017} and \citet{GoekeDis} were the first 
to make note of that QSS reductions are not necessarily the result of a singular 
perturbation in their development of TFPV theory. The insight of 
\citet{GOEKEpara,GOEKE2017} and \citet{GoekeDis} was to recognize 
that if a QSS reduction is the result of a singular perturbation scenario, then 
there \textit{must} exist a normally hyperbolic critical manifold in the singular 
limit. Simply stated, a TFPV is a point in parameter space at which the 
mass action system is equipped with a normally hyperbolic critical manifold 
of non-isolated equilibrium points. 

The IAZA reaction is equipped with a four-dimensional parameter space, 
$\pi := [k_1\;k_{2}\;k_{-1}\;E_T]^T$. Within this parameter space there are 
two TFPVs of significant interest:
\begin{subequations}
\begin{align}
\pi_1^{\star}&:=[0\;k_2\;k_{-1}\;E_T],\\
\pi_2^{\star}&:=[k_1\;0\;k_{-1}\;E_T].
\end{align}
\end{subequations}
The one-dimensional critical manifold associated with $\pi_1^{\star}$ is the 
$z$-axis ($c=0$): $M_0:=\{((z,c)\in \mathbb{R}^2_{\geq 0}:c=0\}$. It is 
straightforward to verify that $M_0$ is normally hyperbolic, and thus the 
decomposition (\ref{split}) holds for all $(z,0)\in M_0 \cap \Lambda$.

To calculate the QSS reduction that corresponds to small $k_1$, we map 
$k_1 \mapsto \varepsilon k_1^{\star}$, and rewrite the mass action equations 
as
\begin{subequations}\label{rescale1}
\begin{align} 
    \dot{z} &= -\varepsilon k_1^{\star}  (E_T - z)z + (k_{-1} + 2\varepsilon k_1^{\star}z)c,\\
    \dot{c} &= \;\;\varepsilon k_1^{\star}(E_T -z)z - (k_{-1}+k_2 + 2\varepsilon k_1^{\star}z)c.
\end{align}
\end{subequations} 
With $0< \varepsilon \ll 1$ and $k_1^{\star}$ carrying unit dimension, we 
express the system (\ref{rescale1}) in singular perturbation form:
\begin{equation}\label{SPF}
    \begin{bmatrix}
    \dot{z} \\ \dot{c}
    \end{bmatrix} = \begin{bmatrix}k_{-1}\\-(k_{-1}+k_2)\end{bmatrix}c + \varepsilon\begin{bmatrix} -k_1^{\star}(E_T-z)z +k_1^{\star}zc\\\;\;k_1^{\star}(E_T-z)z-k_1^{\star}zc\end{bmatrix}.
\end{equation}
By inspection of (\ref{SPF}), we identify the following:
\begin{equation}\label{ingred}
    N:=\begin{bmatrix}k_{-1}\\-(k_{-1}+k_2)\end{bmatrix}, \quad f:=c, \quad Df:=[0\;1], \quad G(z,c,0):=\begin{bmatrix} -k_1^{\star}(E_T-z)z +k_1^{\star}zc\\\;\;k_1^{\star}(E_T-z)z-k_1^{\star}zc\end{bmatrix}.
\end{equation}
From (\ref{ingred}) the projection matrix and associated QSS reduction are 
both straightforward to calculate,
\begin{equation}
    \Pi^{M}:=\begin{bmatrix}1 & \cfrac{K_S}{K_M}\\ 0 & 0\end{bmatrix}, \qquad\begin{bmatrix}\dot{z}\\\dot{c}\end{bmatrix}= \Pi^{M}G(z,0,0)=\begin{bmatrix} -\cfrac{k_1^{\star}k_2}{k_{-1}+k_2}(E_T-z)z\\0\end{bmatrix}.
\end{equation}
Thus, in the limit as $k_1\to 0$, the QSS reduction for $z$ is given by
\begin{equation}\label{ZsQSSA}
    \dot{z} = -\cfrac{k_2}{K_M}(E_T-z)z,
\end{equation}
and from this point forward will refer (\ref{ZsQSSA}) as the standard 
QSS approximation for $z$.

\begin{remark}
Note that the QSS reduction (\ref{ZsQSSA}) is different than the classical 
reduction given by (\ref{zclass}). However, they are asymptotically equal 
in the limit of small $k_1$.
\citet{GOEKE2017} presents a thorough discussion on the asymptotic relationship 
between classical and Fenichel reduction.
\end{remark}

A similar calculation can be performed for $\pi^{\star}_2$ by mapping 
$k_2 \mapsto \varepsilon k_2^{\star}$ and repeating the same procedure. We will 
omit the details and simply state the result (see Appendix for further details). 
In the limit as $k_2\to 0$ with 
all other parameters bounded away from zero, the corresponding QSS reduction is
\begin{equation}\label{pQSSA}
  \dot{z} = -\cfrac{k_2z(E_T-z)(K_S+2z)}{K_S^2+(E_T+2z)K_S+2z^2}, \quad K_S=k_{-1}/k_1,
\end{equation}
and from this point forward we refer to (\ref{pQSSA}) as the QSS approximation 
for $z$ when $w$ is slow. 

Alternatively, in $(w,c)$--coordinates the system
\begin{subequations}
\begin{align}
    \dot{c} &= k_1(e_0+w-c)(z_0-c-w)-(k_{-1}+\varepsilon k_2^{\star})c\\
    \dot{w} &= \varepsilon k_2^{\star} c
\end{align}
\end{subequations}
is in standard form when $k_2$ is small, and the QSS reduction for $w$ is 
\begin{equation}\label{eq:wT}
    \dot{w} = k_2 h^-(w;K_S),
\end{equation}
where ``$h^-(w;K_S)$" in (\ref{eq:wT}) is understood to mean that $K_M$ is 
replaced with $K_S$ in the expression for $h^-(w)$.

What remains is to quantify the approach to the QSS. This can 
be done by analyzing the \textit{layer} problem, obtained by setting 
$\varepsilon =0$ in (\ref{rescale1}). For $k_1=0$, the relationship
\begin{equation}
z-z_0 = \cfrac{K_S}{K_M}(c-c(0))
\end{equation}
holds in the approach to the critical manifold. Thus, if $c(0)=0$, then the 
initial condition for the reduced problem is $z_0$, which agrees with the 
classical reduction. Consequently, the initial value problem
  \begin{equation}\label{ZsQSSA2}
    \dot{z} = -\cfrac{k_2}{K_M}(E_T-z)z, \quad z(0)=z_0,
\end{equation}
is the leading order approximation for $z$ for \textit{small} $k_1$. Furthermore,
integration of (\ref{ZsQSSA2}) yields closed-form solutions for $z$ and $w$:
\begin{subequations}
\begin{align}
z(t) &= \cfrac{z_0E_T}{e_0\exp[{\displaystyle k_2\,E_T t/K_M}]+z_0},\label{eq:sQSSAz}\\
w(t) &= z_0-z(t).\label{eq:sQSSAw}
\end{align}
\end{subequations}

When $k_2=0$, $w$ is conserved in the approach to the critical manifold, 
and $z+c=z_0+c_0$. If $c_0=0$, then $z=z_0-\lambda_Z$ at the onset of the
QSS.\footnote{Here, $\lambda_Z$ is computed with $k_2$ identically $0$.} 
Consequently, the initial value problems
\begin{subequations}\label{pQSSA2}
\begin{align}
  \dot{z} &= -\cfrac{k_2z(E_T-z)(K_S+2z)}{K_S^2+(E_T+2z)K_S+2z^2}, \quad z(0)=z_0-\lambda_Z,\\
  \dot{w} &= k_2h^-(w;K_S), \qquad \qquad \qquad \quad \;\;w(0)=0,
  \end{align}
\end{subequations}
respectively approximate $z$ and $w$ during the QSS phase of the reaction. 
Note that $z(0)$ for (\ref{pQSSA2}) may be much less than $z_0$.

\subsection{The missing piece of the puzzle: \textit{what is $\varepsilon$?}}
The TFPV theory developed by \citet{GOEKEpara,GOEKE2017} provides a direct method 
for obtaining QSS reductions that emerge as a result Fenichel theory and a singular 
perturbation scenario. Hence, it partially answers the origin problem of: 
\textit{where do QSS reductions come from?} However, the entire story in 
incomplete on least three counts. \textit{First}, the notion of ``small" is, at 
this juncture, somewhat abstract and not very quantitative. 
The answer to the question ``How small must $k_1$ be so that (\ref{ZsQSSA}) provides 
a reliable approximation to the full system~(\ref{z5a})-(\ref{z5c})?" is unclear at 
this point, since the term \textit{small} is relative: a concentration is only small 
in comparison to another, larger concentration; a rate constant is \textit{small} only 
in comparison to some other, larger, rate constant. The answer to the question 
pertaining to the magnitude of a perturbation to the TFPV is important, especially 
if the objective 
is to prepare an experiment in a way that ensures (\ref{ZsQSSA}) is accurate. 
Ultimately, if one one wishes to employ (\ref{ZsQSSA}) to estimate kinetic data 
from an assay, then it is necessary that $\pi$ be very close to $\pi_1^{\star}$. 
But, any notion of ``close" remains undefined at this juncture.

\text{Second}, GSPT is not the only method capable of generating reduced models. 
As we already mentioned, QSS parameters defined by \citet{GOEKE2017} establish a 
clear criterion for the near-invariance of the QSS manifold, and near-invariance 
is often sufficient to validate a QSS reduction~\cite{Eilertsen2021}. Hence, not 
all reductions originate from a ``singular perturbation scenario". In addition to 
Fenichel reduction, center manifold reduction is another popular 
technique~\cite{Carr1981,Guckenheimer1983}. The center manifold method is similar 
to Fenichel reduction in that the former reduces the full model via projection 
onto a low-dimensional invariant manifold (the center manifold). However, center manifold 
reduction and Fenichel reduction are not the same: center manifold reduction 
requires the existence of \textit{one} non-hyperbolic equilibrium point, whereas 
Fenichel reduction requires a normally hyperbolic critical manifold, comprised 
of an \textit{infinite} number of equilibria. Since TFPV theory establishes 
criteria for the applicability of Fenichel reduction, it is possible that other 
reductions may exist, and may not originate from a singular perturbation but 
a center manifold scenario (i.e., at a point in parameter space where a 
stationary point in the phase--plane becomes non-hyperbolic). 

\textit{Third}, TFPV theory assumes, by definition, that the critical manifold is 
normally hyperbolic everywhere. This is somewhat restrictive, since in most 
applications one must deal with points where 
\begin{equation}
    \text{rank}\;DfN|_{x\in M_0} < n-k.
\end{equation}

In the model equations of the IAZA reaction mechanism, the loss of normal hyperbolicity 
occurs in the critical set that materializes when
\begin{equation}
    \pi =\pi^{\ddagger}:=[k_1\;0\;0\;E_T].
\end{equation}
In $(w,c)$-space, the critical \textit{set},\footnote{The usage of the word set 
here is intentional, as $\mathcal{S}_1 \cup \mathcal{S}_2$ does not define a manifold.}
$\mathcal{S}_0$, that emerges when $\pi = \pi^{\ddagger}$ is
\begin{equation}\label{tcrit}
    \mathcal{S}_0:= \mathcal{S}_1\cup \mathcal{S}_2=\{(w,c)\in \mathbb{R}^2_{\geq 0}: c=z_0-w\} \cup \{(w,c)\in \mathbb{R}^2_{\geq 0}: c=e_0+w\}.
\end{equation}
If $z_0>e_0$, then the sets $\mathcal{S}_1$ and $\mathcal{S}_2$ intersect in the 
first quadrant at the point $(w,c)=(w_T,E_T/2)$, and the splitting (\ref{split}) 
does not hold. Fortunately, the loss of normal hyperbolicity does not mean we 
cannot construct a QSS reduction; it simply means we must exercise caution in 
the neighborhood of $\mathcal{S}_1 \cap \mathcal{S}_2$. As we will show later 
on, near the singular limit corresponding to $\pi^{\ddagger}$, typical trajectories 
initially follow $\mathcal{S}_2$ but later follow $\mathcal{S}_1$ for the 
remainder of the reaction. Thus, the qualitative behavior of the system is obvious. 
The problem emerges when we try to make quantitative statements about the efficacy 
of the QSS reduction near the singular point at $\mathcal{S}_1 \cap \mathcal{S}_2$. 
It is well-established from Fenichel theory that a small perturbation to the 
vector field $h(x)$ results in the birth of an invariant slow manifold, 
$M_{\varepsilon}$, whose Hausdorff distance, $d_H$, is $\mathcal{O}(\varepsilon)$ 
from $M_0$ (see \cite{kuehn2015multiple}, Chapters 2-3, for a clear presentation 
of general Fenichel theory and its application to fast/slow dynamical systems), 
which is a quantitative statement concerning 
the asymptotic validity of the QSS reduction. Since classical Fenichel theory 
breaks down at $\mathcal{S}_1 \cap \mathcal{S}_2$, the integrity of a QSS 
reduction in the vicinity of this point is unclear.

\section{Unveiling $\varepsilon$: Dimensionless singular perturbation parameters}\label{DTFPV}
In this section we provide an answer to the question: How small must a TFPV in order 
to ensure the associated QSS reduction is reliable? A direct approach involving the
eigenspectrum of championed by \citet{PALSSON1984} ensures the eigenvalue ratio 
of $DH(x)$\footnote{Note that $DH(x)$ is the Jacobian of the full system, whereas 
$Dh(x)$ is the Jacobian of the layer problem.} at the unique equilibrium point is 
small. However, this approach seems rather limited for two reasons. First, the layer 
problem with $\varepsilon =0$ contains an infinite number of equilibrium points, 
and thus estimating the spectral gap at just one of the equilibria may, in some 
cases, be an oversimplification. Second, the eigenspectrum of the Jacobian can 
completely vanish at points where the critical manifold fails to be normally 
hyperbolic, and linear methods do not suffice in such situations. 

The more general approach that we take leverages the phase plane geometry of the 
IAZA reaction. We begin with the preliminary assumption that
\begin{equation}\label{eq:tQSSA}
    c = h^-(w) + \mathcal{E}_Z,
\end{equation}
where $\mathcal{E}_Z$ is an error term. The motivation behind this choice is that 
the QSS manifold approaches the critical manifold(s) in the singular limits discussed 
in Section~\ref{GSPT}. The objective is then to establish criteria 
that ensures $\mathcal{E}_Z$ is small. To generate the condition(s) that render 
$\mathcal{E}_Z$ negligible, we compute the limit supremum
\begin{equation}
    \limsup_{t\to \infty}\; \mathcal{E}_Z^2 \leq \varepsilon^2 \lambda_Z^2.
\end{equation}
The relationship between $e_0$ and $z_0$ alters the value of $\lambda_Z$ 
according to the phase plane geometry described in Section~\ref{sec:phase-plane}. 
As a reminder, $\lambda_Z$ is given by
\begin{equation*}
  \max c := \lambda_Z =  
   \begin{cases}
      h^-((z_0 - e_0) / 2) & \text{if} \quad e_0 < z_0,\\
    h^-(0) & \text{if} \quad  z_0 \leq e_0. 
    \end{cases}
\end{equation*} 
Let us next define
\begin{align*}
    \mu := \sqrt{(K_M + E_T)^2 - 4e_0z_0} \quad \text{ and } \gamma := \sqrt{K_M(K_M + 2E_T)}, 
\end{align*}
so $\lambda_Z$ can be written as
\begin{align*}
 \lambda_Z = \begin{cases} \big(K_M+E_T - \gamma \big)\big/ 2 & \text{if} \quad e_0 < z_0, \vspace{1mm} \\
 \big(K_M+E_T - \mu\big)\big/2 &  \text{if} \quad z_0 \leq e_0. \\
 \end{cases}
 \end{align*} 
The derivative\footnote{The term $D_w$ denotes the derivative with respect to $w$ $D_w:=\text{d}/\text{d}w$.} 
of $\mathcal{E}_Z^2$ is:
\begin{subequations}\label{eq:ENERGY}
\begin{align}
    \cfrac{1}{2}\cfrac{\text{d}\mathcal{E}_Z^2}{\text{d}t} &= [\dot{c} - D_wh^-(w)\dot{w}] \mathcal{E}_Z\\
    &= [k_1(c-h^+(w))(c-h^-(w)) -  D_wh^-(w)\dot{w}]\mathcal{E}_Z.
\end{align}  
\end{subequations}
Phase plane geometry dictates that $c-h^+(w)< 0$ for any trajectory starting 
in or on $\Lambda^*$. It follows from (\ref{eq:ENERGY}) that
\begin{equation}\label{eq:INEQ}
\cfrac{1}{2}\cfrac{\text{d}\mathcal{E}_Z^2}{\text{d}t} \leq  k_1\max[c-h^+(w)]\mathcal{E}_Z^2 + \max|D_wh^-(w)|\max |\dot{w}||\mathcal{E}_Z|.
\end{equation}

\begin{remark}
We have exploited the fact that the derivative function for $\dot{c}$ factors as
$\dot{c}=k_1(c-h^-(w))(c-h^+(w))$ to recover (\ref{eq:INEQ}). It is this geometric 
feature, which is common in most mass actions systems that model enzyme kinetics, 
that allows us to generate sharp error estimates.
\end{remark}

The term ``$\max[c-h^+(w)]$" is easily bounded above by
\begin{equation*}\label{eq:maxs}
\lambda_Z - \displaystyle \min_{0\leq w \leq z_0} h^+(w),
\end{equation*}
where $\min h^+(w)$ is given by
\begin{equation*}\label{eq:maxs2A}
    \displaystyle \min_{0\leq w \leq z_0} h^+(w) = \begin{cases}
      h^+((z_0-e_0)/2) \qquad \text{if} \;\;e_0 < z_0,\\
      h^+(0) \qquad \qquad \qquad \; \text{if} \;\;z_0 \leq e_0.
    \end{cases}
\end{equation*}
Using definition \eqref{eq:maxs2A}, we calculate $\max [c-h^+(w)]$ 
\begin{equation}\label{37}
\max_{c\in \Lambda^{\star}}\; [c-h^+(w)]:=-\varphi=\begin{cases}
-\sqrt{K_M(K_M+2E_T)} = - \gamma \qquad \qquad \;\;\; \text{if} \;\; e_0\leq z_0,\\
\\
-\sqrt{(K_M+E_T)^2-4e_0z_0} =  -\mu \qquad \quad \text{if}\;\; z_0 < e_0.
\end{cases} 
\end{equation}
The last term on the right hand side of (\ref{eq:INEQ}) is easily bounded in terms of 
$\mathcal{E}_Z^2$,
\begin{equation}
  \max |D_wh^-(w)|\max|\dot{w}| |\mathcal{E}_Z| \leq \delta \mathcal{E}_Z^2 +  \cfrac{( \max |D_wh^-(w)| \max|\dot{w}|)^2}{4\delta}.
\end{equation}
Choosing $\delta = k_1\varphi/2$ and applying Gr\"{o}nwall's lemma yields
\begin{subequations}\label{eq:bounds1}
\begin{align}
\mathcal{E}_Z^2 &\leq \mathcal{E}(0)_Z^2 e^{-\displaystyle k_1 \varphi t} + \cfrac{( \max |D_wh^-(w)| \max|\dot{w}|)^2}{(k_1\varphi)^2}\bigg(1-e^{\displaystyle -k_1 \varphi t}\bigg),\\
&\leq \mathcal{E}(0)_Z^2 e^{-\displaystyle k_1 \varphi t} + \cfrac{( \max |D_wh^-(w)| \max|\dot{w}|)^2}{(k_1\varphi)^2}.
\end{align}
\end{subequations}
The constants that bound the $\limsup$ of $\mathcal{E}_Z^2$ are easily computed 
over $\Lambda^*$ using standard calculus techniques
\begin{equation}\label{maxs}
\max_{(c,w)\in \Lambda^*}|D_wh^-(w)| = \cfrac{E_T}{K_M+E_T}, \qquad \max_{(c,w)\in \Lambda^*} |\dot{w}| = k_2 \lambda_Z.
\end{equation}
Finally, from (\ref{maxs}) and (\ref{eq:bounds1}) we have 
\begin{proposition}\label{prop1} 
For every solution of with initial value in $\Lambda^*$ it holds that
\begin{equation}
  \mathcal{E}_Z^2(t) \leq  \mathcal{E}^2(0)e^{-\tau} + \varepsilon^2\lambda_Z^2,
\end{equation}
with $\tau:=k_1\varphi t$. Thus with

\begin{equation}\label{epsL}
    \varepsilon := \bigg(\cfrac{E_T}{K_M+E_T}\bigg)\bigg(\cfrac{K}{K_M}\bigg)
                \bigg(\cfrac{K_M}{\varphi}\bigg), \quad \text{\rm with} \quad K=k_2/k_1 , 
\end{equation}
the solution approaches the QSS manifold, $c=h^-(w)$, up to an error of $\lambda_Z^2 \varepsilon^2$, 
with time constant $\lambda_{\tau}:=k_1\varphi$.
\end{proposition}
\medskip

{\bf Proposition \ref{prop1}} suggests that $\varepsilon$ defined by 
(\ref{epsL}) \textit{is} the singular perturbation parameter that regulates the 
accuracy of the classical reduction~(\ref{wclass}). This follows from the 
observation that a singular perturbation parameter proportional to a ratio of 
fast and slow timescales is \textit{dimensionless}. If one non-dimensionalizes 
the mass action equations, then it is natural to scale $c$ by its maximum
value,\footnote{For excellent references 
on scaling procedures, we invite the reader to 
consult \cite{segel1984modeling,segel1988,shoffner2017approaches}.} 
$\lambda_Z$. In this case the magnitude of the \textit{scaled} error,
$\lambda_Z^{-2}\cdot[c-h^-(w)]^2 := [\hat{c}-\hat{h}^-(w)]^2$, is regulated
(asymptotically) by $\varepsilon$, and therefore one would naturally take 
$\varepsilon$ in {\bf Proposition \ref{prop1}} to be the appropriate singular 
perturbation parameter when the mass action equations are dimensionless. 

The form of $\varepsilon$ in (\ref{epsL}) is revealing, as it factors into
\textit{three} small parameters, $\varepsilon_1$, $\varepsilon_2$ and
$\varepsilon_{\ddagger}$,
\begin{equation}\label{epss}
    \varepsilon_1 := \cfrac{k_1E_T}{k_{-1}+k_2+k_1E_T} \leq 1,\qquad \varepsilon_2:=\cfrac{k_2}{k_{-1}+k_2} \leq 1, \qquad \varepsilon_{\ddagger}:=\cfrac{k_{-1}+k_2}{k_1\varphi} \leq 1.
\end{equation}

\begin{figure}[htb!]
  \centering
 \includegraphics[width=14cm]{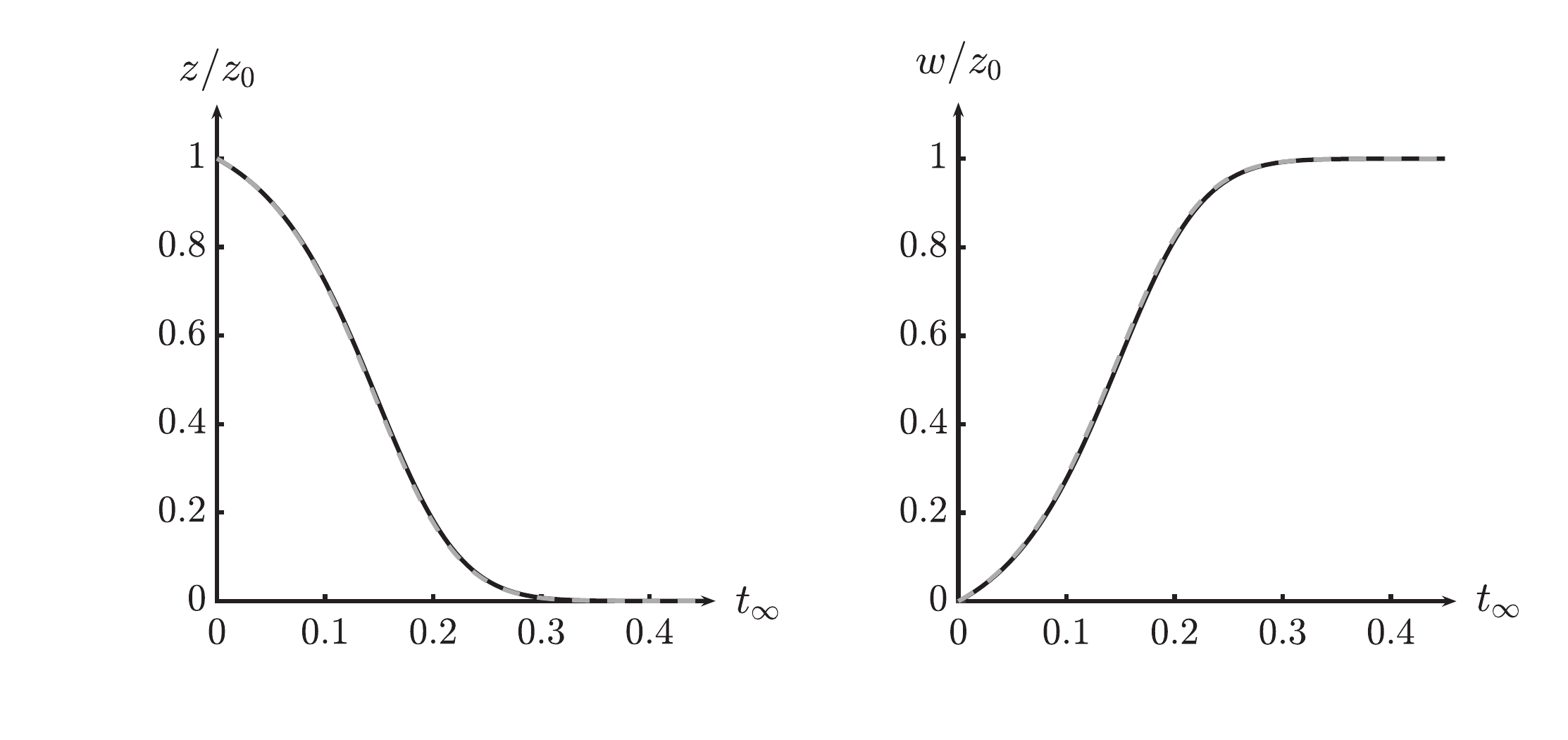}
  \caption{\textbf{The sQSSA, (\ref{eq:sQSSAz}) and (\ref{eq:sQSSAw}), 
  for the IAZA reaction mechanism~(\ref{z1}) is valid when $\mathbf{E_T \ll K_M}$}. In 
  both panels, the solid black curve is the numerical solution to the mass action 
  equations~(\ref{z5a})--(\ref{z5c}), and the dashed black curve (barely visible) 
  is the numerical solution to (\ref{eq:sQSSAz}). The rate constants used 
  to construct the numerical approximation are: $e_0=9$, $z0=1$, $k_1=2$, $k_{-1}=k_2=500$. 
  {\sc Left}: In this panel, the solid black curve is the numerical solution to 
  $z(t)$, and the dashed black curve is the solution generated by (\ref{eq:sQSSAz}). 
  {\sc Right}: In this panel, the solid black curve is the numerical solution 
  to $w(t)$, and the dashed black curve is the solution generated by 
  (\ref{eq:sQSSAw}). For illustrative purposes, the units of all parameters and
  concentrations are arbitrary. Time has been mapped to the $t_{\infty}$ scale: 
  $t_{\infty}(t) =  1-1/\ln[t+\exp(1)]$.}\label{fig44}
\end{figure}

\begin{figure}[htb!]
  \centering
  \includegraphics[width=14cm]{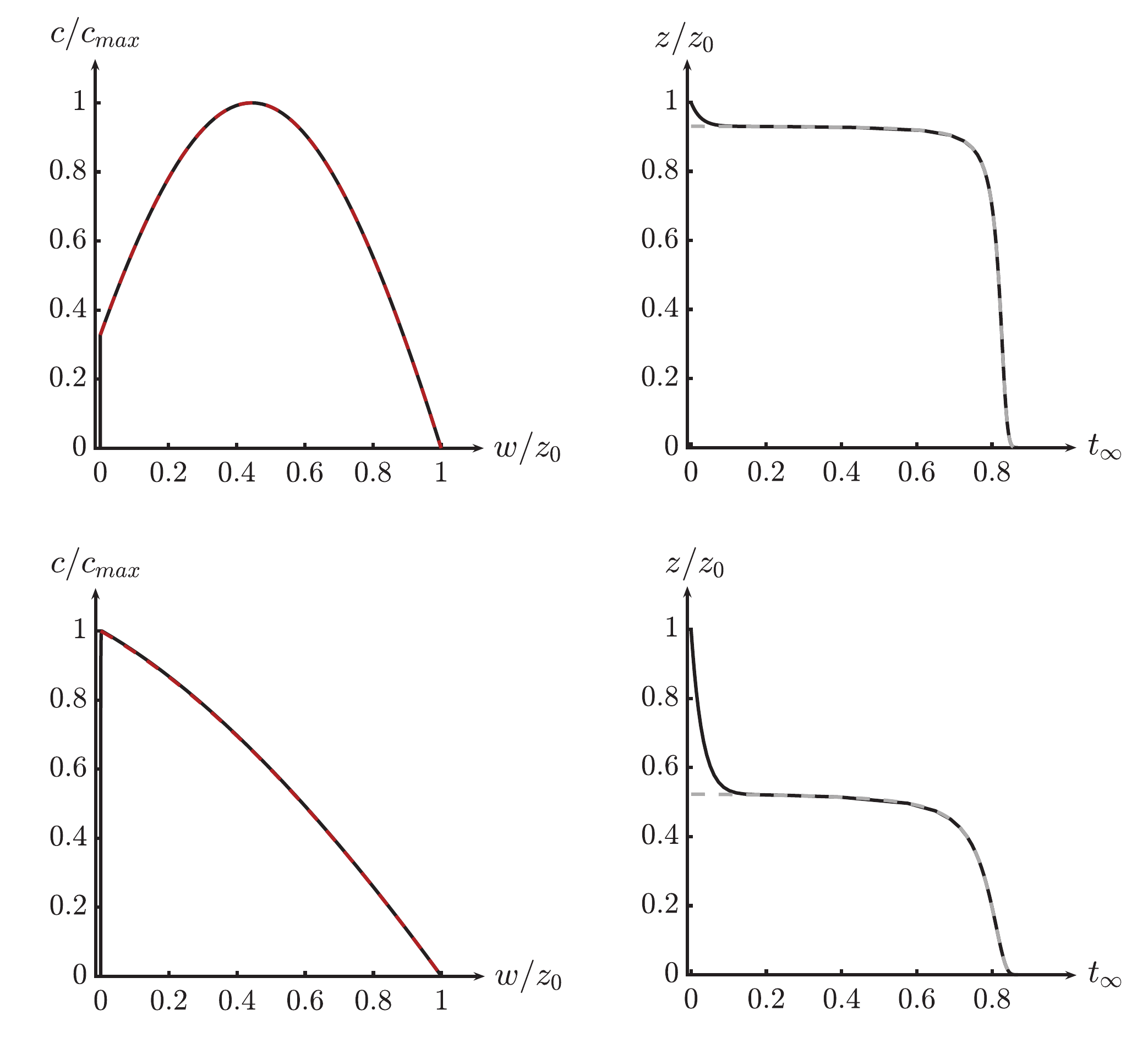}
  \caption{\textbf{The reductions (\ref{eq:wT}) and (\ref{pQSSA}) are valid for 
  the IAZA reaction mechanism~(\ref{z1}) when $\mathbf{k_2 \ll k_{-1}}$.} In all 
  panels, the solid black curve is the numerical solutions to the mass action
  equations~(\ref{z5a})--(\ref{z5c}). In the upper and lower right panels, the 
  dashed/dotted curve is the numerical solution to (\ref{pQSSA}) with $z^*$ 
  as the initial condition. In the top panels, the parameter values and initial
  conditions are $k_1=1.0$, $k_{-1}=5$, $k_2=0.01$ and $e_0=1$, $z_0=9$. In the 
  bottom panels, the parameter values and initial conditions are $k_1=1.0$, 
  $k_{-1}=5$, $k_2=0.01$ and $e_0=7$, $z_0 = 3$. For illustrative purposes, the 
  units of all parameters and concentrations are arbitrary. In the upper and 
  lower right panels, time has been mapped to the $t_{\infty}$-scale: 
  $t_{\infty}(t) =  1-1/\ln[t+\exp(1)]$.
  {\sc Top Left}: This panel illustrates a typical trajectory in   the $(w,c)$ 
  phase plane when the tQSSA is valid. The dashed red curve is the $c$-nullcline $h^-(w)$, 
  and by inspection it is easy to see that (\ref{eq:wT}) is justified when 
  $k_2 \ll k_{-1}$. {\sc Top Right}: In this panel, the numerical solutions to 
  $z$ as well as (\ref{pQSSA}) are compared on the $t_{\infty}$ scale. Clearly,
  (\ref{pQSSA}) serves as a good approximation to $z(t)$. The lower left and 
  lower right panels give analogous illustrations of the $(w,c)$ phase plane and
  temporal dynamics when $e_0\geq z_0$.}\label{fig5}
\end{figure}

\begin{remark}
By inspection, it is clear that conditions $\varepsilon_1 \ll 1$ and 
$\varepsilon_2 \ll 1$ determine the validity of the standard QSS approximation 
for $z$ (\ref{ZsQSSA}) and the QSS approximation for $z$ when $w$ is slow
(\ref{pQSSA}), respectively. This follows from the observation that 
$\varepsilon_1=0$ if $k_1=0$, and $\varepsilon_2=0$ if $k_2=0$, and thus the 
limiting cases of $\varepsilon_1\to 0$ and $\varepsilon_2 \to 0$ correspond 
with the TFPVs $\pi_1^{\star}$ and $\pi_2^{\star}$. Thus, the 
\textit{dimensionless parameters} $\varepsilon_1$ and $\varepsilon_2$ give us 
a measure of the ``size" of $k_1$ or $k_2$ with respect to a perturbed 
$\pi_1^{\star}$ or $\pi_2^{\star}$. Numerical results support this claim; 
see {{\sc Figures}~\ref{fig44} and \ref{fig5}}.
\end{remark}

At this juncture, it is important to give the dimensionless parameters $\varepsilon_1$ 
and $\varepsilon_2$ a proper interpretation in the context of singular perturbation 
reduction and TFPV theory. Recall that the reduction (\ref{ZsQSSA}) is valid --- and 
justified from singular perturbation theory --- as $k_1 \to 0$. The question we set 
out to answer pertained to how small one must take $k_1$ in order to be confident 
that (\ref{ZsQSSA}) accurately approximates the full system on the slow timescale. 
Let us now rephrase this question with more clarity: if all other parameters are 
fixed, or bounded \textit{below} and \textit{above} by positive constants, how 
small, in comparison to the other parameters, must we take $k_1$ to secure an 
accurate reduction in (\ref{ZsQSSA})? The answer is until 
$\varepsilon_1 \ll 1$. If $E_T$, $k_{-1}$ and $k_2$ are held constant, then 
$k_1$ should be reduced in magnitude until $k_1$ is at least small enough that 
$\varepsilon_1 \ll 1$. Loosely speaking, if one were to define a one-dimensional 
search direction in parameter space along $k_1$, with all other parameters 
constant, then $\varepsilon_1 \ll 1$ can be interpreted as a stopping criterion. 
A similar interpretation holds for $k_2$, $\varepsilon_2$, and (\ref{pQSSA}) 
and (\ref{pQSSA2}). 

Note that $\varepsilon_1 \to 0$ vanishes not only as $k_1 \to 0$, but also as 
$E_T\to 0$. However, the latter does not result in a singular perturbation 
scenario,\footnote{With $E_T$=0, the vector field is void of a critical manifold; 
hence Fenichel theory does not apply.} and thus one cannot attribute the accuracy 
of (\ref{zclass}) to Fenichel theory. Nevertheless, {\bf Proposition~\ref{prop1}} 
clearly indicates that (\ref{ZsQSSA}) will improve as $E_T\to 0$, and this demands 
explanation. First, note that the normal form of a \textit{static} transcritical 
bifurcation
\begin{equation}\label{TBnormalform}
    \dot{x} = rx-x^2,
\end{equation}
is recovered from (\ref{ZsQSSA}) by reparameterizing time, $\tau := -k_2t/K_M$, 
and thus $E_T$ in (\ref{ZsQSSA}) is equivalent to the bifurcation parameter ``$r$" 
in the normal form (\ref{TBnormalform}). Moreover, setting $E_T=0$ in the rate 
equations results in the formation of a non-hyperbolic fixed point at the origin. 
Consequently, the origin is equipped with center manifold, $W^c(0)$, when $E_T=0$. 
For \textit{small} $E_T$ and $z\in[0,z_0]$, the long-time dynamics can be 
approximated via projection onto the local center manifold, $W^c_{\text{loc.}}(0)$, 
defined in the \textit{extended} phase space, $(z,c,E_T)\in \mathbb{R}^3.$ To 
leading order, the dynamics of $z$ on the center manifold is
\begin{equation*}
    \dot{z}= -\cfrac{k_2}{K_M}(E_T-z)z,
\end{equation*}
which is equal to (\ref{ZsQSSA}) (see Appendix for the detailed calculation). This 
example illustrates that the mathematical justification of the QSS reduction may 
be multifaceted. Ultimately, the justification depends on which \textit{path} is 
taken in parameter space. For example, taking $k_1 \to 0$ with all other parameters 
constant (and positive) results in justification from Fenichel theory. Center manifold 
reduction legitimizes the QSS reduction as $E_T\to 0$ with all other parameters 
positive and constant and $z \leq z_0 < E_T$. Thus, while the reductions as 
$k_1\to 0$ and $E_T \to 0$ are the same, the \textit{origin} is different.

The remaining dimensionless parameter, $\varepsilon_{\ddagger}$, vanishes when 
both $k_{-1}$ \textit{and} $k_2$ vanish as long as $z_0 < e_0$, which is consistent with $\pi^{\ddagger}$. 
However, $\pi^{\ddagger}$ is not a TFPV due 
to the lack of normal hyperbolicity. Nevertheless --- and on an intuitive level ---
we expect $\varepsilon_{\ddagger}$ to tell us something about the efficacy of 
a QSS reduction near $\pi^{\ddagger}$. Such is the subject of 
Section~\ref{sec:rQSSA}.

\section{The loss of normal hyperbolicity}\label{sec:rQSSA}
\subsection{The Layer Problem}
In perturbation form, and in $(c,w)$ coordinates, the mass action system for 
small $k_{-1}$ and $k_2$ is
\begin{equation}\label{reverse}
\begin{bmatrix}
    \dot{c} \\ \dot{w}
    \end{bmatrix} = \begin{bmatrix} k_1(e_0+w-c)(z_0-c-w) \\ 0\end{bmatrix} + \varepsilon  \begin{bmatrix}-(k_{-1}^{\star}+k_2^{\star})\\k_2^{\star}\end{bmatrix}c,
\end{equation}
and is thus a fast/slow system in standard form. The layer problem corresponds 
to the singular limit with $\varepsilon =0$ in (\ref{reverse}), but the critical 
set (\ref{tcrit}) of the layer problem is not a manifold. Moreover, normal
hyperbolicity does not hold at $\mathcal{S}_1 \cap \mathcal{S}_2$ located 
at $(w,c) = ((z_0-e_0)/2,E_T/2):= (w_T,c_T)$, since 
\begin{equation}
\cfrac{\partial}{\partial c}\;k_1(e_0+w-c)(z_0-c-w)\bigg|_{(c,w)=(c_T,w_T)}=0.
\end{equation}
However, apart from $\mathcal{S}_1\cap \mathcal{S}_2$, the stability of the 
additional equilibria that comprise $\mathcal{S}_0$
is straightforward to compute
\begin{subequations}\label{stability}
\begin{align}
\cfrac{\partial}{\partial c}\;k_1(e_0+w-c)(z_0-c-w)\bigg|_{c=z_0-w} &= -k_1(e_0-z_0+2w):=\lambda_1(w),\\
\cfrac{\partial}{\partial c}\;k_1(e_0+w-c)(z_0-c-w)\bigg|_{c=e_0+w} &= \;\;k_1(e_0-z_0+2w):=\lambda_2(w),
\end{align}
\end{subequations} 
and from (\ref{stability}) a clear picture emerges:

\begin{proposition}\label{prop2}
With $\mathcal{S}_1 :=\{(c,w)\in \mathbb{R}^2_{\geq 0} : c=z_0-w\}$,
$\mathcal{S}_2 :\{(c,w)\in \mathbb{R}^2_{\geq 0}: c=e_0+w\}$, and $\delta_1,\delta_2,\delta_3,\delta_4>0$, 
the following holds:\footnote{The use of $\delta_1$,...,$\delta_4$ allow us to define 
$\mathcal{S}_{1,2}^{a,r}$ as \textit{compact} submanifolds.}
\begin{enumerate}
\item If $z_0 > e_0$ and $w > w_T$, then $\lambda_2(w) >0$ and 
$\mathcal{S}_2^r:=\{(c,w) \in \mathbb{R}_{\geq 0}^2 \cap \mathcal{S}_2 : w \geq w_T +\delta_1 \}$ 
is comprised of repulsive fixed points w.r.t the layer problem given by the 
singular limit of (\ref{reverse}).

\item If $z_0 > e_0$ and $w > w_T$, then $\lambda_1(w) <0$ and 
$\mathcal{S}_1^a:=\{(c,w) \in \mathbb{R}_{\geq 0}^2 \cap \mathcal{S}_1 : w \geq w_T+\delta_2 \}$ 
is comprised of attracting fixed points with respect to the layer problem 
given singular limit of (\ref{reverse}).

\item If $z_0 > e_0$ and $w < w_T$, then $\lambda_2(w) <0$ and 
$\mathcal{S}_2^a:=\{(c,w) \in \mathbb{R}_{\geq 0}^2 \cap \mathcal{S}_2 : w \leq w_T-\delta_3 \}$ 
is comprised of attracting fixed points with respect to the layer problem 
given by the singular limit of (\ref{reverse}).

\item If $z_0 > e_0$ and $w < w_T$, then $\lambda_1(w) > 0$ and 
$\mathcal{S}_1^r:=\{(c,w) \in \mathbb{R}_{\geq 0}^2 \cap \mathcal{S}_1 : w \leq w_T -\delta_4 \}$ 
is comprised of repulsive fixed points with respect to the layer problem 
given by the singular limit of (\ref{reverse}).
\end{enumerate}
\end{proposition}

\begin{figure}[htb!]
  \centering
    \includegraphics[width=9cm]{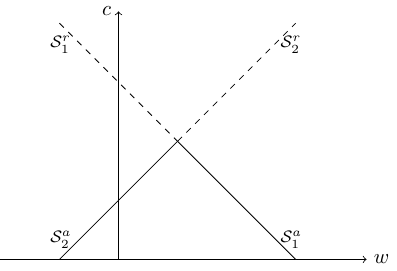}
  \caption{\textbf{In the singular limit that occurs when $\mathbf{K_M=0}$, 
  a dynamic transcritical bifurcation emerges in the vector field}. When
  $k_{-1}=k_2=0$, the nullclines $c=h^-(w)$ and $ c=h^+(w)$ coalesce with 
  the curves $\mathcal{S}_1:=\{(c,w)\in \mathbb{R}^2_{\geq 0}:c=z_0-w\}$ and
  $\mathcal{S}_2 :=\{(c,w)\in \mathbb{R}^2_{\geq 0} :c=w+e_0\}$. Each curve, 
  $c=z_0-w$ and $c=w+e_0$, is comprised of stationary points. These
  curves intersect and exchange stability at the transcritical bifurcation 
  point $(w,c) = ((z_0-e_0)/2,E_T/2) \equiv(w_T,c_T)$. Thus, $\mathcal{S}_1$ 
  and $\mathcal{S}_2$ have stable (attracting) and unstable (repulsive) regions,
  $\mathcal{S}_1^a,\mathcal{S}_2^a$ and $\mathcal{S}_1^r,\mathcal{S}_2^r$,
  respectively. Dashed lines correspond to unstable equilibrium points, and 
  solid lines correspond to stable equilibrium points. }\label{fig6}
\end{figure}

\begin{figure}[htb!]
  \centering
    \includegraphics[width=8.0cm]{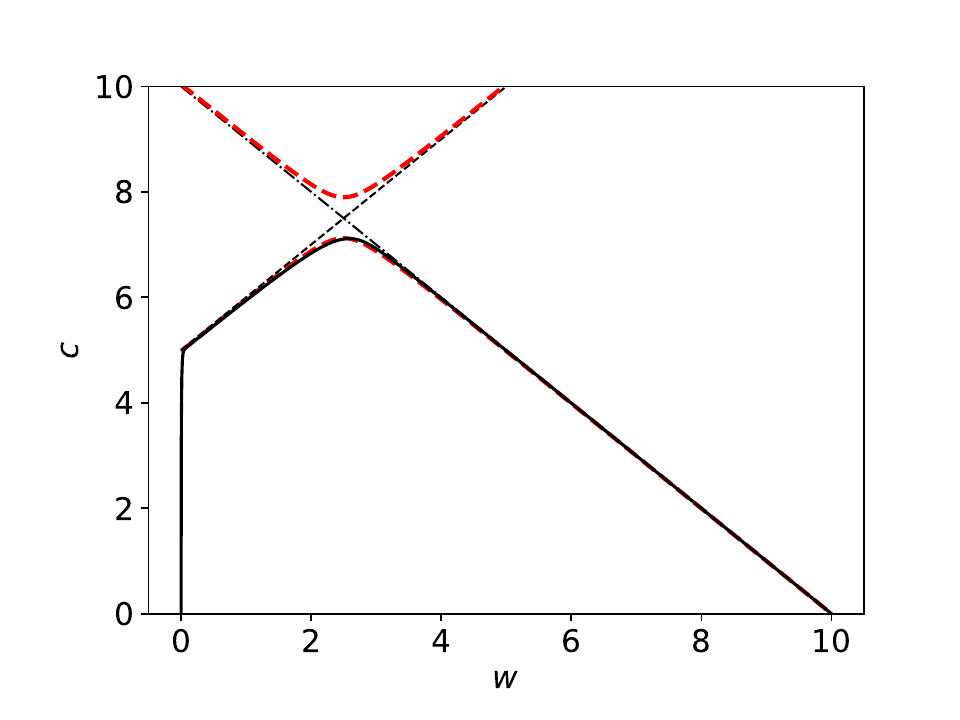}
    \includegraphics[width=8.0cm]{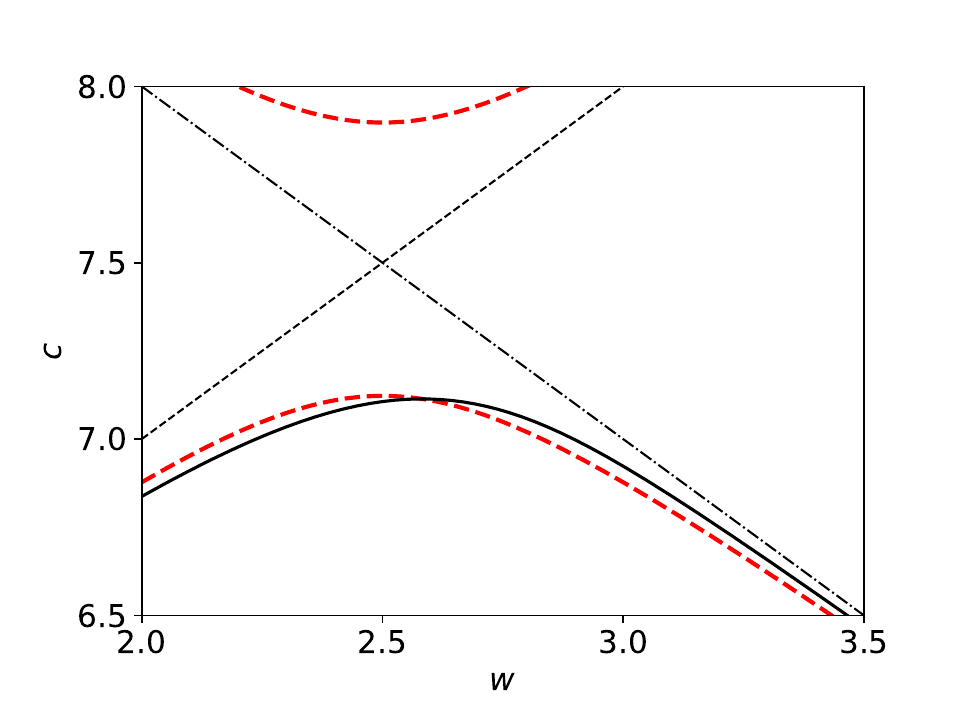}
  \caption{\textbf{Fenichel theory breaks down at $\mathcal{S}_1 \cap \mathcal{S}_2$
  and the accuracy of reduced models in the neighborhood of this point must be
  estimated}. In both panels, the thick black curve is the numerical solution to 
  the mass action system (\ref{z5a})--(\ref{z5c}). The thick dashed red curves are
  $h^{\pm}(w)$. The thin dashed-dotted black line is $\mathcal{S}_1$, and the 
  thin, dashed black line is $\mathcal{S}_2$. Parameters values (of arbitrary units) 
  are: $e_0=5$, $z_0=10$, $k_{-1}=0.1$, $k_2=0.1$, and $k_1 =10.0$.  
 } \label{rQSSApanels}
\end{figure}

Thus, from {\bf Proposition \ref{prop2}} the geometric landscape of the layer problem 
is clear: the two submanifolds $\mathcal{S}_1$ and $\mathcal{S}_2$ exchange 
stability at their intersection (see {\sc Figure}~\ref{fig6}).
The question that follows is: It is possible to construct a useful QSS reduction 
and, if so, how do we quantify its efficacy? If $e_0 < z_0$, then physical 
solutions that originate in $\Lambda^*$ will follow $\mathcal{S}_1^a$ and
$\mathcal{S}_2^a$, and the difficulty lies in gauging the accuracy of either 
the Fenichel reduction or the classical reduction (\ref{wclass}) near the point 
of intersection, $\mathcal{S}_1\cap \mathcal{S}_2$ (see 
{{\sc Figure}}~\ref{rQSSApanels}).
One could resort to center manifold reduction but, the eigenspectrum at 
$\mathcal{S}_1 \cap \mathcal{S}_2$ contains a zero eigenvalue with multiplicity 
two and, as noted in \cite{jardonkojakhmetov2019survey}, center manifold reduction 
is not very useful in such scenarios. However, from a 
practical point of view, note that if $z_0 < e_0$, then the layer problem has 
only one \textit{attracting} branch ($\mathcal{S}_1$) contained entirely within 
the first quadrant. Thus, with $z_0 < e_0$, the intersection 
$\mathcal{S}_1 \cap \mathcal{S}_2$ lies in the second quadrant, and the reduced 
problem is
\begin{equation}\label{rQSSA}
\dot{w}=k_2(z_0-w).
\end{equation}

In enzyme kinetics, the linear QSS reduction (\ref{rQSSA}) is the autocatalytic 
version of the \textit{reverse} quasi-steady-state approximation 
(rQSSA)~\cite{Segel1989,2000SchnellMaini,EILERTSEN2020}. The experimental utility of 
the rQSSA is that the catalytic rate constant, $k_2$, can be estimated by fitting 
(through regression analysis) the exact solution of (\ref{rQSSA}) to time course 
data~\cite{choiandKim2017}. Of course, one will normally want to refine experiments 
by taking different $z_0$ with fixed $e_0$ and compute an average over the range 
of $k_2$ estimates. Thus, to design a proper experimental protocol capable of estimating
$k_2$, we must answer the twofold question: How small should $k_{-1}$ and $k_2$ 
be and, practically speaking, what ratios $z_0/e_0$ should be considered? As 
we illustrate in the subsection that follows, the answers to these questions 
are inexorably linked to Fenichel theory.

\subsection{Investigating normal hyperbolicity in the irreversible IAZA reaction}
To gain insight pertaining to the validity of (\ref{rQSSA}), we will temporarily 
study the completely irreversible IAZA reaction, which is analogous to the 
Van Slyke-Cullen reaction mechanism:
\begin{equation}\label{z1R}
    \ce{$Z$ + $E$ ->[$k_1$] $C$ ->[$k_2$] $2E$ + $W$}.
\end{equation}
In the singular limit that coincides with $k_2=0$, the geometry of the layer problem 
of (\ref{z1R}) is unchanged (see {\sc Figure}~\ref{fig6}). Once the perturbation is
``turned on" and $0 < k_2$, a distinct invariant trajectory, $\mathcal{C}^{\varepsilon}$,
emerges
\begin{equation}
    \mathcal{C}^{\varepsilon}:=\mathcal{S}_1.
\end{equation}
Thus, the set $\mathcal{S}_1$ remains invariant. Moreover, since the trajectory
$\mathcal{C}^{\varepsilon}$ travels from the unstable ($\mathcal{S}_1^r)$ 
to the stable ($\mathcal{S}_1^a$) sections of $\mathcal{S}_1$, one can refer to this 
distinguished trajectory as a \textit{faux Canard.} 

There are two cases we must consider in our assessment of reduced models obtained 
via projection onto normally hyperbolic and attracting subsets of 
$\mathcal{S}$: $z_0\leq e_0$ and $e_0 <z_0$. We will begin with the latter. 
Linearizing about $x^{(1)}$ with $e_0=z_0$ yields
\begin{equation*}
\begin{bmatrix} \dot{c}\\\dot{w}\end{bmatrix} =DH(c,w)|_{(0,z_0)}=\begin{bmatrix}-2e_0k_1-k_2 & -2e_0k_1 \\ k_2 & 0\end{bmatrix}\begin{bmatrix} c\\w\end{bmatrix}.
\end{equation*}
The eigenvalues of $DH(0,z_0)$ are $\lambda_1:=-k_2$ and $\lambda_2:=-2k_1e_0$, and 
the corresponding eigenvectors are
\begin{equation*}
    v_1,v_2:=\begin{bmatrix} -1\\\;\;\;1 \end{bmatrix}, \quad \begin{bmatrix} -2e_0k_1/k_2\\ 1 \end{bmatrix}.
\end{equation*}
If $k_2\ll 2e_0k_1$, then $\text{span} \;v_1$ is the slow eigenspace. In the slow 
timescale, $T=k_2t$, the linear solution is
\begin{equation*}
    x=\beta_1 v_1 e^{-\displaystyle T} + \beta_2 v_2 e^{-\displaystyle T/\varepsilon}, \qquad \varepsilon:= \lambda_1/\lambda_2, \qquad x:=[c \;\; w-z_0]^T,
\end{equation*}
where $\beta_1$, $\beta_2$ are constants. For initial conditions that are sufficiently 
close to $x^{(1)}$, the solution will lie almost entirely along the direction of 
$v_1$ after a time of order $T=-\varepsilon \ln \varepsilon$.

Certainly $\lambda_1/\lambda_2\ll 1$ is a necessary condition for the validity of 
(\ref{rQSSA}). However, since the linearized solution is only valid in a small 
neighborhood surrounding $x^{(1)}$, it is unclear at this point if (\ref{rQSSA}) 
holds after a timescale of order $T=-\varepsilon\ln \varepsilon$ for initial 
conditions that lie sufficiently far away from $x^{(1)}$. For a nonlinear fast/slow 
system of the form
\begin{subequations}\label{FastSLow}
\begin{align}
    x' &= f(x,y),\\
    \varepsilon y' &= g(x,y),
\end{align}
\end{subequations}
with $\mathfrak{Re}(D_yg(x,y))<0 \;\;\forall (x,y)\in M_0$, trajectories will reach 
an $\mathcal{O}(\varepsilon)$--neighborhood of the slow manifold once 
$T\sim -\varepsilon \ln \varepsilon$. However, due to the lack of normal 
hyperbolicity at $\mathcal{S}_1\cap \mathcal{S}_2$, a timescale estimate of 
$T\sim -\varepsilon\ln \varepsilon$ is possibly too short when $e_0=z_0$. 
Consequently, it is
necessary to get a more global estimate on how quickly trajectories approach 
$\mathcal{C}^{\varepsilon}$ when $z_0 = e_0$ if (\ref{rQSSA}) is to be useful.
\begin{proposition}\label{prop3}
For all $(c,w)\in \Lambda$ with $z_0\leq e_0$ and $\varepsilon_{\ddagger}>0$, it holds that
\begin{equation*}
\mathcal{Z}^2\leq \mathcal{Z}^2(0) e^{\displaystyle (\varepsilon_{\ddagger}-1)T/\varepsilon_{\ddagger}},
\end{equation*}
where $\mathcal{Z}^2:=[c-(z_0-w)]^2$ and $T=k_2t$.
\end{proposition}

\noindent \textit{Proof}: For $z_0\leq e_0$, we have
\begin{subequations}
\begin{align*}
\cfrac{\text{d}\mathcal{Z}^2}{\text{d}t} &\leq -2k_1(e_0-\lambda_Z)\mathcal{Z}^2,\\
&=-k_1(e_0-z_0+\sqrt{K^2+2E_TK}-K)\mathcal{Z}^2,\\
&\leq-k_2K^{-1}(\sqrt{K^2+2E_TK}-K)\mathcal{Z}^2.
\end{align*}
\end{subequations}
As defined in the previous section, $K:=k_2/k_1$. Passing to the slow time, $T$, we obtain:
\begin{subequations}
\begin{align*}
 \cfrac{\text{d}\mathcal{Z}^2}{\text{d}T}&\leq   -K^{-1}(\sqrt{K^2+2E_TK}-K)\mathcal{Z}^2,\quad T:=k_2t,\\
 &=(1-\varepsilon_{\ddagger}^{-1})\mathcal{Z}^2,\\
 &= (\varepsilon_{\ddagger}-1)\varepsilon_{\ddagger}^{-1}\mathcal{Z}^2,
\end{align*}
\end{subequations}
and Gr\"{o}nwall's lemma implies
\begin{equation*}
    \mathcal{Z}^2 \leq \mathcal{Z}^2(0) e^{\displaystyle (\varepsilon_{\ddagger}-1)T/\varepsilon_{\ddagger}} \approx \mathcal{Z}^2(0) e^{\displaystyle -T/\varepsilon_{\ddagger}}\quad \text{for}\quad \varepsilon_{\ddagger}\ll 1.
\end{equation*}

\begin{remark}
It is important to realize that $\varepsilon_{\ddagger}$ is in fact undefined when $e_0 = z_0$ and $k_{-1} = k_2 = 0$. This is a result of the loss of normal hyperbolicity. However, with $k_{-1} = 0$ and $k_2 \to 0$, the limit from the right does exist (it is zero), and this is all that matters in applications since we are bounded away from the singular limit.
\end{remark}

{\bf Proposition~\ref{prop3}} confirms our intuition from Section~\ref{DTFPV} that
$\varepsilon_{\ddagger}$ is the singular perturbation parameter that regulates 
the accuracy of the rQSSA when $z_0 \leq e_0$. Moreover, it is important to emphasize 
that $\varepsilon_{\ddagger}$ is different, in an important way, from the singular 
perturbation parameter 
that emerges from scaling analysis. To see this, let $e_0=z_0$ and define
$\bar{c}:=c/e_0$ and $\bar{w}:=w/e_0$. The dimensionless mass action system is
\begin{subequations}\label{dimson}
\begin{align}\label{dim}
    \varepsilon^{\star}\bar{c}' &= (1+\bar{w}-\bar{c})(1-\bar{c}-\bar{w})-\varepsilon^{\star}\bar{c},\\
    \bar{w}' &= \bar{c},
\end{align}
\end{subequations}
where ``$\phantom{x}'$" denotes differentiation with respect to the slow time, $T$, 
and $\varepsilon^{\star}:=k_2/k_1e_0$. However, $\varepsilon^{\star}$ is not 
asymptotically equivalent to $\varepsilon_{\ddagger}$ when $z_0=e_0$. In fact, 
with $e_0=z_0$, and $E_T=2e_0$ we have
\begin{equation}\label{ORDERS}
    \varepsilon_{\ddagger}\approx \cfrac{1}{2}\sqrt{\varepsilon^{\star}}\sim \mathcal{O}(\sqrt{\varepsilon^{\star}}) \quad \text{for} \quad k_2\ll k_1e_0.
\end{equation}

\begin{figure}[htb!]
  \centering
    \includegraphics[width=8.0cm]{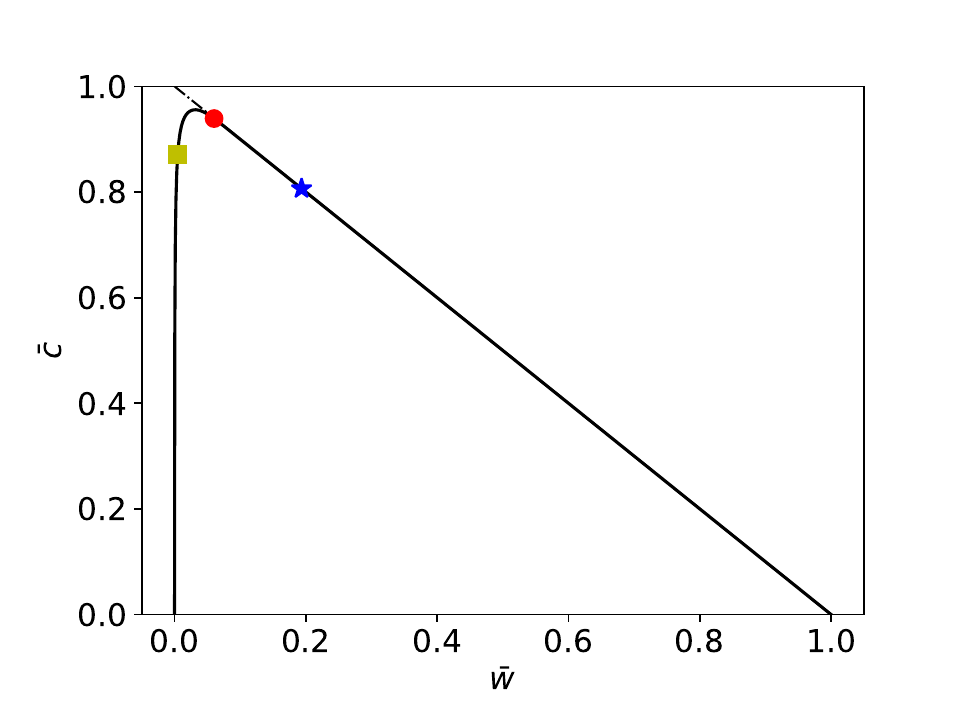}
    \includegraphics[width=8.0cm]{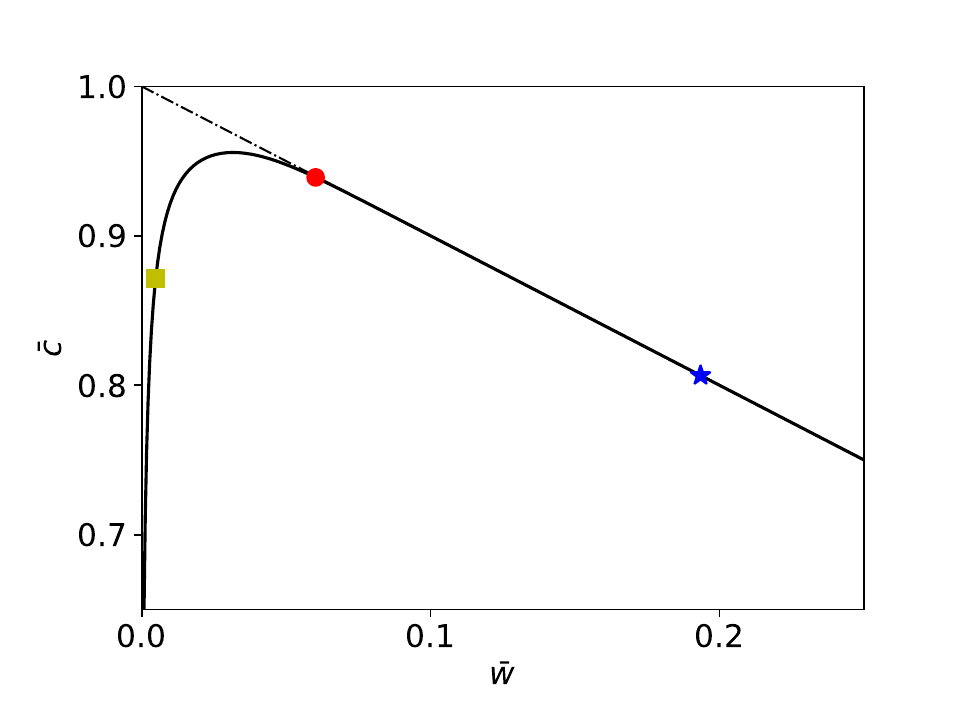}
  \caption{\textbf{Near the bifurcation point, it takes longer than 
  $T\sim -\varepsilon^{\star}\ln \varepsilon^{\star}$ to reach an 
  $\mathcal{O}(\varepsilon^{\star})$--neighborhood of the slow manifold.} In both 
  panels, the thick black curve is the numerical solution to the dimensionless 
  system~(\ref{dimson}) with $\varepsilon^{\star}=0.001$ and 
  $\bar{c}(0)=\bar{w}(0)=0$. The thin, dash-dotted line is $\mathcal{S}_1$. The 
  yellow square demarcates that position of the trajectory when 
  $T=-\varepsilon^{\star}\ln \varepsilon^{\star}$. The red circle shows when 
  $T= -\varepsilon_{\ddagger}\ln \varepsilon_{\ddagger}$ and the blue star at
  $T=-\sqrt{\varepsilon^{\star}}\ln \varepsilon^{\star}$. By inspection, it is 
  clear that at $T=-\varepsilon^{\star}\ln \varepsilon^{\star}$, the trajectory 
  is still not very close to $\mathcal{S}_1$. The {{\sc right}} panel is a close-up 
  of the {{\sc left}}   panel near the bifurcation point.} \label{times}
\end{figure}

The relationship between $\varepsilon^{\star}$ and $\varepsilon_{\ddagger}$ given 
by (\ref{ORDERS}) allows us to better estimate the time it takes for a phase plane 
trajectory that starts at an $\mathcal{O}(1)$ distance from the slow manifold 
to reach an $\mathcal{O}(\varepsilon^{\star})$--neighborhood of the slow manifold:
\begin{equation}\label{sqr}
    |\bar{c}-(1-\bar{w})| \approx e^{\displaystyle -T/2\varepsilon_{\ddagger}} \approx e^{\displaystyle -T/\sqrt{\varepsilon^{\star}}}.
\end{equation}
We see from (\ref{sqr}) that the dimensionless system (\ref{dim}) may require 
a time that is possibly closer to the order 
$T\sim -\sqrt{\varepsilon^{\star}}\ln \varepsilon^{\star}$ before the trajectory is
$\mathcal{O}(\varepsilon^{\star})$ from the slow manifold. In addition, 
{\bf Proposition~\ref{prop3}} indicates that after a time of roughly,
\begin{equation}
    T \sim -\varepsilon_{\ddagger}\ln \varepsilon_{\ddagger},
\end{equation}
the trajectory should be within an 
$\mathcal{O}(\varepsilon_{\ddagger})\sim \mathcal{O}(\sqrt{\varepsilon^{\star}})$--neighborhood 
of the slow manifold. Thus, we \textit{estimate}
\begin{equation}
    -\varepsilon^{\star}\ln \varepsilon^{\star} < T^{\varepsilon^{\star}}\leq -\sqrt{\varepsilon^{\star}}\ln\varepsilon^{\star} \quad \text{for} \quad z_0= e_0\quad \text{and}\;\; k_2\ll k_1e_0,
\end{equation}
where $T^{\varepsilon^{\star}}$ is dimensionless the time it takes before the 
trajectory is within an $\mathcal{O}(\varepsilon^{\star})$--neighborhood from 
the slow manifold. Numerical simulations confirm this estimate (see 
{{\sc Figure}}~\ref{times}). 

Now for the second scenario with $e_0 < z_0$. Defining $\widehat{c}:=2c/E_T$, we 
recover the \textit{general} dimensionless system:
\begin{subequations}\label{scaledG}
\begin{align}
    \widehat{\varepsilon}\;\widehat{c}'&=(\sigma + \theta \bar{w}-\widehat{c})(\theta -\widehat{c}-\theta\bar{w})-\widehat{\varepsilon}\;\widehat{c},\\
    \bar{w}'&=\theta^{-1}\widehat{c}, 
\end{align}
\end{subequations}
where $\widehat{\varepsilon}:=2k_2/k_1E_T$, $\sigma:=2e_0/E_T$ and $\theta:=2z_0/E_T$. 
If $\varepsilon_{\ddagger}\ll 1$ and $e_0 < z_0$, the phase plane trajectory will 
rapidly approach $\mathcal{S}_2^a$, pass near the bifurcation point, and then follow
$\mathcal{S}_1^a$ for the remainder of the reaction. Thus, the accuracy of the reduced 
problem is still limited to its accuracy in the neighborhood of the bifurcation point.

\begin{proposition}
For any trajectory starting in $\Lambda^*$, the asymptotic validity of the reduced 
problem corresponding to (\ref{scaledG}) is limited by a term of 
$\mathcal{O}(\sqrt{\widehat \varepsilon}\;)$ in the neighborhood of the bifurcation 
point. 
\end{proposition}

\noindent \textit{Proof}:
The distance, $d_{\mathcal{B}}$, from $\mathcal{S}_1 \cap \mathcal{S}_2$ to $h^-(w_T)$ 
is identically $\sqrt{K^2+2E_TK}-K$. Every trajectory starting in $\Lambda^*$ is 
bounded away from the bifurcation point by $d_{\mathcal{B}}$, which is bounded 
above by $\sqrt{2E_TK}$:
\begin{equation*}
   d_{\mathcal{B}}= \sqrt{K^2+2E_TK}-K \leq \sqrt{2E_TK} \sim \sqrt{K}.
\end{equation*}
Scaling this bound by $E_T/2$ (the maximal possible distance between $h^-(w_T)$ 
and $\mathcal{S}_1 \cap \mathcal{S}_2$) yields
\begin{equation}\label{square}
   \cfrac{2\sqrt{2E_TK}}{E_T} = 2\sqrt{\cfrac{2K}{E_T}} \sim \mathcal{O}(\sqrt{\widehat{\varepsilon}}\;).
\end{equation}

Our estimate (\ref{square}) indicates that the scaling law\footnote{By ``scaling law", 
we are referring to the asymptotic order of which the trajectory tracks the slow 
manifold.} of the dimensionless system (\ref{scaledG}) is proportional to the 
square root of the eigenvalue ratio is in agreement with established theoretical 
results. Trajectories of fast/slow systems in standard form~(\ref{FastSLow}) pass 
near dynamic transcritical bifurcation points at a distance that is of order
$\sqrt{\varepsilon}$ (see \cite{Berglund} for a detailed overview of scaling laws 
near bifurcation points). The blow-up 
method \cite{Dumortier1993,Krupa2001,kuehn2015multiple} could possibly be utilized 
to rigorously explain the dynamics of the rQSSA (\ref{rQSSA}) and the classical 
reduction (\ref{wclass}) in neighborhoods containing 
$\mathcal{S}_1\cap \mathcal{S}_2$. However, such an analysis is beyond the scope 
of the paper, and will be the subject of a future work.

One practical challenge we face, both computationally and experimentally, is how 
to properly handle the \textit{switching} undergone by the trajectory as it 
moves from one attracting submanifold to the other. To circumnavigate this 
difficulty, we utilize the classical reduction (\ref{wclass}). From 
{\bf Proposition~\ref{prop1}} we have
\begin{equation*}
    \varepsilon :=\bigg(\cfrac{k_1E_T}{k_2 + k_1E_T}\bigg)\cdot \bigg(\cfrac{K}{\sqrt{K(K+2E_T)}}\bigg),
\end{equation*}
which holds for all $e_0 \leq z_0$. The \textit{upper bound}
\begin{subequations}\label{UB}
\begin{align*}
   [\hat{c}-\hat{h}^-({w})]^2 &\leq [\hat{c}-\hat{h}^-({w})]^2(0)e^{\displaystyle -k_1\varphi t} + \varepsilon^2,\\
   &=[\hat{c}-\hat{h}^-({w})]^2(0)e^{\displaystyle -k_1\varphi\cdot \cfrac{k_2}{k_2}\cdot t} + \varepsilon^2,\\
   &=[\hat{c}-\hat{h}^-({w})]^2(0)e^{\displaystyle - T/\varepsilon } + \varepsilon^2
\end{align*}
\end{subequations}
follows from {\bf Proposition \ref{prop1}}, and with $k_2 \ll k_1E_T$, $\varepsilon$ of 
{\bf Proposition~\ref{prop1}} is approximately
\begin{equation}\label{asym}
    \varepsilon \approx  \cfrac{K}{\sqrt{K(K+2E_T)}} \approx \sqrt{K/2E_T}\quad \text{for small} \;\;K.
\end{equation}
We conclude from (\ref{asym}) that the classical reduction (\ref{wclass}) is 
valid as long as $\sqrt{K/2E_T}\ll 1$. 

\subsection{A geometric interpretation of reversibility and the validity of the rQSSA}
\begin{figure}[htb!]
  \centering
    \includegraphics[width=10cm]{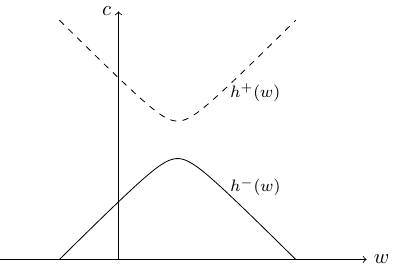}
  \caption{\textbf{Reversibility may be interpreted as an imperfection}. Adding 
  reversibility by setting $k_{-1}>0$ in the layer problem with $k_{2}=0$ yields 
  two disjoint critical manifolds, $c=h^-(w)$ and $c=h^+(w)$. As indicated by the 
  solid line, the critical manifold $c=h^-(w)$ is comprised of attracting fixed 
  points. In contrast, the critical manifold $c=h^+(w)$ (dashed line) is comprised 
  of repulsive stationary points. Note the bifurcation structure observed in 
  {{\sc Figure}} \ref{fig6} is broken unless $k_{-1}=0$.} \label{fig7}
\end{figure}

When reversibility is added to the binding step of the reaction mechanism~(\ref{z1R})
so that $k_{-1}\neq 0$ but $k_2=0$ in the singular limit, the layer problem consists 
of two disjoint critical manifolds, $c=h^+(w)$ and $c=h^-(w)$, and the set
$\mathcal{C}^{\varepsilon}$ is no longer invariant. In this way, reversibility has 
an interpretation as a dynamic imperfection, since the bifurcation structure reported 
in {{\sc figure} \ref{fig6}} is broken (see {{\sc figure}} \ref{fig7}). 

From {\bf Proposition~\ref{prop2}}, we see that the classical reduction's accuracy 
improves with increasing $k_{-1}$. Moreover, the rQSSA (\ref{rQSSA}) can still be
successfully employed when $z_0\leq e_0$ provided $k_{-1}\ll k_1e_0$, and the 
validity of the rQSSA when $z_0=e_0$ is still regulated by a term that is 
proportional to $\sqrt{K_M/E_T}$. Again, the geometric interpretation is 
that $d_{\mathcal{B}}$ is
\begin{equation}
   d_{\mathcal{B}} =\sqrt{K_M^2+2E_TK_M} -K_M
\end{equation}
when $0<k_{-1}$, and thus the \textit{scaled} distance is proportional to $\sqrt{K_M/E_T}$, 
while the total distance is ultimately proportional to $\sqrt{2E_TK_M}$. Since the 
distance from $h^-(w)$ to $\mathcal{S}_1$ decreases as $w \to z_0$ (i.e., because 
the stationary point $(c,w)=(0,z_0)$ lies in the graphs of both $\mathcal{S}_1$ 
and $h^-(w)$), we expect the rQSSA (\ref{rQSSA}) to improve as the ratio of $z_0$ 
to $e_0$ diminishes. Again, the blow-up method should be employed if one wishes to 
make more rigorous statements concerning the accuracy of the classical 
reduction and the rQSSA near the bifurcation point when $K_M \ll E_T.$

\begin{remark}
To minimize the influence of the bifurcation point, an experiment should be prepared so 
that $z_0\ll e_0$ and $K_M\ll e_0$. Thus, the experimental qualifier, 
$K_M\ll E_T, z_0\ll e_0$, is more restrictive than the mathematical qualifier, 
$K_M\ll E_T, z_0\leq e_0$. Moreover, if in addition to $K_M\ll E_T$ it holds that 
$k_2\ll k_{-1}$, then the classical reduction (\ref{wclass}) will yield the most 
accurate result; see {\bf Proposition \ref{prop1}}.
\end{remark}

\section{Discussion}\label{Diss}
In this paper, we have found that the mass action equations of the IAZA reaction mechanism 
admit three possible QSS reductions: the standard QSSA~(\ref{ZsQSSA}), the reverse 
QSSA~(\ref{rQSSA}), and ``total" QSSA~(\ref{eq:wT}). Furthermore, we have systematically 
derived the qualifiers (timescale ratios) that warrant the validity of a particular 
reduction (see, {\sc Table~\ref{table:1}}).
\begin{table}[htb!]
\centering
\begin{tabular}{{c}|{c}|{c}|{c}}
\textbf{Parameter}(s) &\textbf{Qualifier} & \textbf{Approximation} &  \textbf{Critical Set} \\
\hline
\hline
$k_1$&\;\;\quad $\varepsilon_1:=k_1E_T/(k_{-1}+k_2) \leq tol.$ & sQSSA (\ref{ZsQSSA}) & $c_1=0$  \\
$k_2,k_{-1}$&\;\;\quad $\widetilde \varepsilon \;\;:= \sqrt{(k_{-1}+k_2)/k_1E_T}\leq tol.,\;\; \&\;\;z_0\leq e_0$           
&  rQSSA (\ref{rQSSA}) & $\mathcal{S}_0:=\mathcal{S}_1\cup \mathcal{S}_2$ \\
$k_2$&\;\;\quad $\varepsilon_2 := k_{2}/k_{-1} \leq tol.$          & tQSSA (\ref{eq:wT}) & $h^-(w;K_S)$ \\
\hline
\end{tabular}
\center
\caption{{\bf The parameters, qualifiers, and critical sets for the standard, 
reverse, and total QSSAs}. The first column lists the dimensional parameters (rate 
constants) that, as a consequence of Fenichel theory, ensure the accuracy of the 
corresponding QSS reduction (third column) as they approach zero with all other 
parameters bounded away from zero. Of course, traditional Fenichel reduction does 
not apply to $\mathcal{S}_0$, but one can compute reductions that are valid along 
compact submanifolds of $\mathcal{S}_0$. The second column provides a metric (or 
stopping condition along the search direction defined by the parameter in the first 
column) that indicates how small one must make the parameter in the first column 
before the associated reduction (third column) is sufficiently accurate. The 
tolerance ``$tol.$" in the second column can be used to define was is meant by 
sufficiently accurate. Mathematically, it is desirable that $tol. \ll 1$. However, 
there is limited resolution in laboratory measurements and, for experiments, 
a $tol.$ of $0.01$ is usually sufficient for an accurate QSS
reduction~\cite{Schnell2003Estimates}. The fourth column lists the critical 
set recovered in the singular limit that emerges when the parameter(s) in the 
first column is identically zero (with all other parameters bounded away 
from zero). }\label{table:1}
\end{table}

On the mathematical side, we have demonstrated that finding ``$\varepsilon$," can 
be a difficult task, as even planar systems can exhibit complicated phenomena, 
such as dynamic bifurcations where the classical results of GSPT are inapplicable. 
Additionally, the mathematical mechanism that warrants a QSS reduction may depend 
on the path one takes in parameter space. Here we have shown that the sQSSA is 
justifiable either from Fenichel reduction (as $k_1 \to 0^+$), or from center 
manifold reduction (as $E_T\to 0^+$), and the mathematical justification depends 
on which parameter vanishes in the limit. 

Our analysis of the IAZA reaction also alludes to a rather interesting research 
topic, which is the possible extension of TFPV theory to include dynamic bifurcation 
points. For example, it is understood that, due to the loss of normal hyperbolicity, 
asymptotic scaling laws change near bifurcation points~\cite{Berglund}. The 
results of our analysis are consistent with the established scaling law near a 
transcritical bifurcation, as we demonstrated that the accuracy of the rQSSA 
when $z_0=e_0$ is proportional to $\sqrt{k_2}$ in the absence of reversibility. 
This suggests that if we define ``$\varepsilon$" to be synonymous with $k_2$ 
via the mapping $k_2 \mapsto \varepsilon k_2^*$, then the accuracy of the rQSSA 
for the irreversible mechanism should scale as $\sqrt{\varepsilon k_2^*}$ near 
$\mathcal{S}_1 \cap \mathcal{S}_2$. While further investigation of this topic 
may appear moot or biochemically irrelevant with respect to experimental assay 
design, it should be noted that many in vivo reactions occur under conditions 
requiring stochastic fluctuations in the model equations. 
In such cases, it is necessary to resort to a Langevin equation or to the 
chemical master equation. Furthermore, QSS reductions of stochastic models 
are often adopted from the non-elementary rate equations generated from the 
deterministic QSS reduction and, interestingly, are frequently accurate under 
the same\footnote{We are abusing the word ``same" here; stochastic rate constants 
differ from deterministic rate constants, and in the chemical master equation
regime the state space 
is discrete.} conditions that validate the deterministic QSS 
reduction~\cite{Sanft,Kang2019,MacNamara2008}. QSS reductions of stochastic models 
are extremely useful in computational biology and computational medicine, as it 
is generally less expensive to simulate the reduced model~\cite{TURNER2004}. 
Consequently, knowledge of both singular perturbation parameters, as well as 
scaling laws near bifurcation points, may help elucidate the validity of QSS 
approximations to stochastic models. 

From the experimental perspective, we note that timescale separation is generally
\textit{induced} in laboratory experiments by preparing the reaction so that the ratio 
of initial enzyme concentration to reactant concentration is small (i.e., with 
$e_0 \ll z_0$); see~\cite{fuentes2005kinetics,Garcia-Moreno1991,Wu2001}. However, 
the condition $e_0\ll z_0$ does not ensure that the spectral 
gap of the Jacobian is wide enough to validate (\ref{zclass}) or (\ref{wclass}). This 
result differs from non-autocatalytic Michaelis Menten-type reactions, where a 
disparity in initial enzyme and substrate concentration is some times enough to ensure 
the validity of the standard QSS reduction~\cite{Heineken1967,Segel1989}. We also 
note that the IAZA reaction mechanism analyzed in the work was comparatively simple, 
and more complicated mechanisms have been proposed~\cite{fuentes2005kinetics}. 
Moreover, most reactions can be controlled with 
inhibitors or activators. In the  the presence of these modifiers, the determination 
of singular perturbation parameters for the reaction mechanism is much more difficult. 
This is because \textit{experimentally practical} QSS reductions of three- and 
four-dimensional systems will generally result from projecting onto a one-dimensional 
slow manifold.\footnote{Generally, time course data is extracted from a single chemical 
species in a laboratory assay, which makes one-dimensional QSS reductions particularly
favorable.} Thus, while it may be that only one eigenvalue vanishes in the singular 
limit, the accuracy of the QSS reduction can require more than one dimensionless 
singular perturbation parameter to be small since there can be multiple fast timescales 
in higher dimensional systems; see \cite{eilertsen2019characteristic} as an example.
Understanding the qualifiers that ensure the accuracy of QSS reductions in 
higher-dimensional systems is critical if we are to make progress in metabolic engineering
and drug development. And, as we have illustrated here, \textit{there is room for 
mathematics.}

\section*{Acknowledgement}
Justin Eilertsen and Santiago Schnell wish to the thank the participants 
of {\em Mathematics of Chemical Reaction Networks}, a seminar jointly 
hosted by the University of Vienna, Politecnico di Torino, and the University 
of Copenhagen, for helpful discussions provided after the presentation of this
work. This manuscript has significantly improved thanks to the insight and careful
comments made by Reviewer~1. Justin Eilertsen's work was partially supported 
by the University of Michigan Postdoctoral Pediatric Endocrinology and Diabetes 
Training Program ``Developmental Origins of Metabolic Disorder'' (NIH/NIDDK 
grant: T32 DK071212). \\

\section*{Appendix}

\subsection{Fenichel reduction in the small $k_2$ limit}
In this appendix, we provide the full calculation for (\ref{pQSSA}). To start, we 
define $k_2 \mapsto \varepsilon k_2^{\star}$, and express the mass-action equations 
in perturbation form:
\begin{equation*}
    \begin{bmatrix}
    \dot{z} \\ \dot{c}
    \end{bmatrix} = \begin{bmatrix} -1 \\ \;\;\;1\end{bmatrix} \bigg(k_1(E_T-z)z+(k_{-1}+2k_1z)c \bigg)+ \varepsilon \begin{bmatrix}0\\k_2^{\star}c\end{bmatrix}.
\end{equation*}
By inspection, we have
\begin{equation*}
    N(z,c):=\begin{bmatrix} -1 \\ \;\;\;1\end{bmatrix}, \quad f(z,c):=k_1(E_T-z)z+(k_{-1}+2k_1z)c, \quad \text{and} \quad G(z,c):=\begin{bmatrix}0\\k_2^{\star}c\end{bmatrix}.
\end{equation*}
The derivative, $Df$ is the row vector
\begin{equation*}
    Df:=[\partial_z f \;\; \partial_cf],
\end{equation*}
where ``$\partial_x$" denotes $\partial/\partial x$. Next, the term $N Df$ is 
an inner product and hence scalar-valued:
\begin{equation*}
    Df N := \langle Df, N\rangle =-\cfrac{\partial}{\partial z}f(z,c) + \cfrac{\partial}{\partial c}f(z,c).
\end{equation*}
The projection operator, $\Pi^M$, is thus
\begin{subequations}
\begin{align*}
    \Pi^M&= \mathbb{I}^{2 \times 2} - N(DfN)^{-1}Df,\\
    &= \mathbb{I}^{2 \times 2} - ((DfN)^{-1})NDf,\\
    \end{align*}
\end{subequations}
and the QSS reduction (\ref{pQSSA}) is
\begin{equation*}
    \dot{x}=\Pi^M G(z,c), \quad \text{evaluated at} \quad c=\cfrac{(E_T-z)z}{(2z+K_S)},
\end{equation*}
with $x$ denoting $[z \; c]^T.$

\subsection{Center manifold reduction in the small $E_T$ limit}
Here we demonstrate that the QSS reduction, (\ref{ZsQSSA}), is justified via 
center manifold reduction as $E_T \to 0^+$. First, we write the extended system 
\begin{equation*}
\begin{bmatrix}\dot{E}_T \\ \dot{z} \\ \dot{c} \end{bmatrix} = \begin{bmatrix}0 & 0 & 0\\ 0 & 0 & k_{-1}\\ 0 & 0 & -(k_{-1}+k_2)\end{bmatrix}\begin{bmatrix}E_T \\ z \\ c\end{bmatrix} + \begin{bmatrix}\;\;\;0 \\ \;\;\;1\\ -1\end{bmatrix}(2k_1zc - k_1(E_T-z)z),
\end{equation*}
which close enough to standard form. At lowest order, the local center manifold 
has the representation
\begin{equation*}
    c=h(z,E_T) = \alpha z^2 + \beta E_Tz + \gamma E_T^2 + \mathcal{O}(z^3,E_Tz^2,E_T^2z,E_T^3).
\end{equation*}
To determine the coefficients $\alpha,\beta$ and $\gamma$, we substitute 
$h(z,E_T)$ into the invariance equation,
\begin{equation*}
    D_z h(z,E_T) \dot{z}(z,h(z,E_T)) + D_{E_T}h(z,E_T)\dot{E}_T = \dot{c}(z,h(z,E)),
\end{equation*}
and equate identical polynomial orders to zero. This yields 
$\alpha = -K_M^{-1}, \beta = K_M^{-1}$ and $\gamma=0.$ Thus, $W_{\text{loc.}}^c(0)$ is given by 
\begin{equation*}
    c\sim h(z,E_T) =\cfrac{k_1(E_T-z)z}{k_{-1}+k_2},
\end{equation*}
and the dynamics of $z$ near the origin (and for small $E_T$) is given by
\begin{equation*}
    \dot{z} = -\cfrac{k_2}{K_M}(E_T-z)z + \mathcal{O}(z^3,E_Tz^2).
\end{equation*}
Consequently, we not only recover (\ref{ZsQSSA}) in the small $E_T$ limit, we see 
that the vector field undergoes a transcritical bifurcation as $E_T$ moves from 
positive to negative.


\end{document}